\documentclass[12pt]{amsart}
\usepackage{amscd,amssymb,epsfig}
\calclayout
\newcommand\tr{\operatorname{tr}}
\newcommand\skw{\operatorname{skw}}

\newcommand\sym{\operatorname{sym}}
\newcommand\grad{\operatorname{grad}}
\renewcommand\div{\operatorname{div}}
\newcommand\curl{\operatorname{curl}}

\newcommand\K{\mathbb{K}}
\newcommand\R{\mathbb{R}}
\newcommand\eps{\operatorname\epsilon}
\newcommand\x{\times}
\newcommand\A{{\mathcal A}}
\renewcommand\P{{\mathcal P}}

\newcommand\J{{\mathcal J}}
\newcommand\E{{\mathcal E}}

\newcommand\Th{{\mathcal T}_h}
\newcommand\T{{\mathbb T}}
\newcommand\M{{\mathbb M}}
\newcommand\V{{\mathbb V}}
\newcommand\W{{\mathbb W}}
\newcommand\X{{\mathbb X}}
\newcommand\Sym{{\mathbb S}}
\newcommand\Hdiv{H(\div,\Omega;\Sym)}
\newcommand\Hdivn{H(\div,\Omega;\M)}

\DeclareMathOperator{\vect}{\operatorname{vec}}

\newcommand\norm[1]{\lVert#1\rVert}

\numberwithin{equation}{section}
 
\theoremstyle{plain}
\newtheorem{thm}{Theorem}[section]
\newtheorem{prop}[thm]{Proposition}
\newtheorem{lem}[thm]{Lemma}

\theoremstyle{remark}
\newtheorem*{remark}{Remark}

\begin{document}
\title[Mixed methods for elasticity with weakly
imposed symmetry]{Mixed finite element methods for linear elasticity 
with weakly imposed symmetry}

\author{Douglas N. Arnold, Richard S. Falk, and Ragnar Winther}

\thanks{October 31, 2005. Revised September 1, 2006.}

\thanks{The work of the first author was supported in part by NSF grant
DMS-0411388.  The work of the second author was supported in part by NSF grant
DMS03-08347. The work of the third author was supported by the Norwegian
Research Council.}

\address{Institute for Mathematics and its Applications, University of
  Minnesota, Minnesota
55455}

\email{arnold@ima.umn.edu}

\urladdr{http://www.ima.umn.edu/\char'176arnold/}

\address{Department of Mathematics, Rutgers University, NJ 08854-8019}

\email{falk@math.rutgers.edu}

\urladdr{http://www.math.rutgers.edu/\char'176falk/}

\address{Centre of Mathematics for Applications and Department of Informatics, 
University of Oslo, P.O. Box 1053, Blindern, 0316 Oslo, Norway}

\email{ragnar.winther@cma.uio.no}

\urladdr{http://folk.uio.no/\char'176rwinther/}



%\renewcommand{\subjclassname}{\textup{2000} Mathematics Subject 
%Classification}
%\subjclass{Primary: 65N30, Secondary: 74S05}

\begin{abstract} 
In this paper, we construct new finite element methods for the
approximation of the equations of linear elasticity in three space
dimensions that produce direct approximations to both stresses and
displacements. The methods are based on a modified form of the
Hellinger--Reissner variational principle that only weakly imposes the
symmetry condition on the stresses. Although this approach has been
previously used by a number of authors, a key new ingredient here is
a constructive derivation of the elasticity complex starting
from the de~Rham complex. By mimicking this construction in the
discrete case, we derive new mixed finite elements for elasticity in a
systematic manner from known discretizations of the de Rham complex.
These elements appear to be simpler than the ones previously derived. 
For example, we construct stable discretizations which use only
piecewise linear elements to approximate the stress field and piecewise
constant functions to approximate the displacement field.
\end{abstract}
\maketitle
\section{Introduction}\label{sec:intro}
The equations of linear elasticity can be written as a system of
equations of the form
\begin{equation}\label{mixed-system}
A\sigma = \eps u, \qquad \div \sigma = f \quad \text{in }\Omega.
\end{equation}
Here the unknowns $\sigma$ and $u$ denote the stress and displacement
fields engendered by a body force $f$ acting on a linearly elastic
body which occupies a region $\Omega \subset \R^3$. Then $\sigma$
takes values in the space $\Sym:=\R^{3\x3}_{\text{sym}}$ of symmetric
matrices and $u$ takes values in $\V:=\R^3$. The differential operator
$\eps$ is the symmetric part of the gradient, the div operator is
applied row-wise to a matrix, and the compliance tensor
$A=A(x):\Sym\to\Sym$ is a bounded and symmetric, uniformly positive
definite operator reflecting the properties of the body.  If the body
is clamped on the boundary $\partial\Omega$ of $\Omega$, then the
proper boundary condition for the system \eqref{mixed-system} is $u =
0$ on $\partial\Omega$. For simplicity, this boundary condition will
be assumed here.  The
issues that arise when other boundary conditions are assumed (e.g.,
the case of pure traction boundary conditions $\sigma n = g$)
are discussed in \cite{acta}.

The pair $(\sigma,u)$ can alternatively be characterized as
the unique critical point of the Hellinger--Reissner functional
\begin{equation}\label{HR}
\J(\tau,v)=\int_\Omega\bigl(\frac12
A\tau:\tau+\div\tau\cdot v-f\cdot v\bigr)\,dx.
\end{equation}
The critical point is sought 
among all $\tau\in\Hdiv$,
the space of square-integrable symmetric matrix fields with
square-integrable divergence, and all $v\in L^2(\Omega;\V)$, the
space of square-integrable vector fields.
Equivalently, $(\sigma,u)\in \Hdiv \x L^2(\Omega;\V)$
is the unique solution to the following weak formulation
of the system \eqref{mixed-system}:
\begin{equation}\label{mixed-system_w}
\begin{array}{lrll}
\int_\Omega(A\sigma:\tau + \div\tau\cdot u)\,dx &=& 0,
\quad
&\tau \in \Hdiv,\\[0.2cm]
\int_\Omega \div\sigma\cdot v \,dx  &=& 
\int_\Omega f \cdot v\,dx, \quad
&v \in L^2(\Omega;\V).
\end{array}
\end{equation}

A mixed finite element method determines an approximate stress field
$\sigma_h$ and an approximate displacement field $u_h$ as the critical
point of $\J$ over $\Sigma_h\x V_h$ where $\Sigma_h\subset\Hdiv$ and
$V_h\subset L^2(\Omega;\V)$ are suitable piecewise polynomial
subspaces. Equivalently, the pair $(\sigma_h,u_h)\in\Sigma_h\x V_h$ is
determined by the weak formulation \eqref{mixed-system_w}, with the test space
restricted to $\Sigma_h\x V_h$.  As is well known, the subspaces
$\Sigma_h$ and $V_h$ cannot be chosen arbitrarily. To ensure that a
unique critical point exists and that it provides a good approximation of the
true solution, they must satisfy the stability conditions from Brezzi's
theory of mixed methods \cite{Brezzi,Brezzi-Fortin}.

Despite four decades of effort, no stable simple mixed finite
element spaces for elasticity have been constructed. For the
corresponding problem in two space dimensions, stable finite elements
were presented in \cite{Arnold-Winther02}.  For the lowest order
element, the space $\Sigma_h$ is composed of piecewise cubic functions,
with 24 degrees of freedom per triangle, while the space $V_h$ consists
of piecewise linear functions. Another approach which has been
discussed in the two dimensional case is the use of composite elements,
in which $V_h$ consists of piecewise polynomials with respect to one
triangulation of the domain, while $\Sigma_h$ consists of piecewise
polynomials with respect to a different, more refined, triangulation
\cite{Arnold-Douglas-Gupta,Fraejisdv,Johnson-Mercier,Watwood-Hartz}. In
three dimensions, a partial analogue of the element in
\cite{Arnold-Winther02} has been proposed and shown to be stable in
\cite{Adams-Cockburn}.  This element uses piecewise quartic stresses
with 162 degrees of freedom per tetrahedron, and piecewise linear
displacements.

Because of the lack of suitable mixed elasticity elements, several authors
have resorted to the use of Lagrangian functionals which are modifications of
the Hellinger--Reissner functional given above
\cite{Amara-Thomas,Arnold-Brezzi-Douglas,Arnold-Falk,Stein-Rolfes,
Stenberg86,Stenberg88,Stenberg-mafelap}, in which the symmetry of the stress
tensor is enforced only weakly or abandoned altogether.  In order to discuss
these methods, we consider the compliance tensor
$A(x)$  as a symmetric and positive definite operator mapping $\M$
into $\M$, where $\M$ is the space of $3 \x 3$ matrices.
In the isotropic case, for example, the mapping $\sigma \mapsto A \sigma$
has the form
\[
A \sigma = \frac{1}{2\mu}\bigl(\sigma - 
\frac{\lambda}{2\mu + 3 \lambda} \tr(\sigma)I\bigr),
\]
where $\lambda(x), \mu(x)$ are positive scalar coefficients, the Lam\'{e}
coefficients.
A modification of the variational principle discussed above is
obtained if we consider the extended Hellinger--Reissner functional
\begin{equation}\label{HRe}
\J_e(\tau,v,q)= \J(\tau,v) + 
\int_\Omega
\tau: q\,dx
\end{equation}
over the space $\Hdivn \x L^2(\Omega;\V)\x L^2(\Omega;\K)$, where $\K$
denotes the space of skew symmetric matrices.  We note that the symmetry 
condition
for the space of matrix fields is now enforced through the introduction of the
Lagrange multiplier, $q$.  A critical point $(\sigma,u,p)$ of the functional
$\J_e$ is characterized as the unique solution of the system
\begin{equation}\label{emixed-system}
\begin{array}{lrll}
\int_\Omega(A\sigma:\tau + \div\tau\cdot u + \tau : p)\,dx &=& 0,
\quad
&\tau \in \Hdivn,\\[0.2cm]
\int_\Omega \div\sigma\cdot v \,dx  &=& 
\int_\Omega f \cdot v\,dx, \quad
&v \in L^2(\Omega;\V),\\[0.2cm]
\int_\Omega \sigma :q \,dx  &=& 0, \quad
&q \in L^2(\Omega;\K).
\end{array}
\end{equation}
Clearly, if $(\sigma,u,p)$ is a solution of this
system, 
then $\sigma$ is symmetric, i.e., $\sigma \in \Hdiv$,  and therefore the pair 
$(\sigma,u) \in \Hdiv \x L^2(\Omega;\V)$ solves the corresponding
system \eqref{mixed-system_w}.  On the other hand, if $(u,p)$ solves
\eqref{mixed-system_w}, then $u\in H^1(\Omega;\V)$ and, if we set
$p$ to the skew-symmetric part of $\grad u$, then $(\sigma,u,p)$ solves
\eqref{emixed-system}.  In this respect, the two systems
\eqref{mixed-system_w} and \eqref{emixed-system} are equivalent.
However, the extended system \eqref{emixed-system} leads to new
possibilities
for discretization. Assume that we choose finite element spaces
$\Sigma_h\x V_h \x Q_h \subset \Hdivn \x L^2(\Omega;\V)\x
L^2(\Omega;\K)$
and consider a discrete system corresponding to
\eqref{emixed-system}.
If $(\sigma_h, u_h,p_h) \in \Sigma_h\x V_h \x Q_h$ is a discrete
solution, 
then $\sigma_h$
will not necessary inherit the symmetry property of $\sigma$.
Instead, $\sigma_h$ will satisfy the weak symmetry condition
\[
\int_\Omega \sigma_h :q \,dx  = 0, \qquad \text{for all } q \in Q_h.
\]
Therefore, these solutions will in general not correspond to 
solutions of the discrete system obtained from \eqref{mixed-system_w}.

Discretizations based on the system \eqref{emixed-system} will be referred to
as mixed finite element methods with weakly imposed symmetry.  Such
discretizations were already introduced by Fraejis de Veubeke in
\cite{Fraejisdv} and further developed in \cite{Arnold-Brezzi-Douglas}.  In
particular, the so--called PEERS element proposed in
\cite{Arnold-Brezzi-Douglas} for the corresponding problem in two space
dimensions used a combination of piecewise linear functions and cubic bubble
functions, with respect to a triangulation of the domain, to approximate the
stress $\sigma$, piecewise constants to approximate the displacements, and
continuous piecewise linear functions to approximate the Lagrange multiplier
$p$. Prior to the PEERS paper, Amara and Thomas \cite{Amara-Thomas} developed
methods with weakly imposed symmetry using a dual hybrid approach.  The lowest
order method they discussed approximates the stresses with quadratic
polynomials plus bubble functions and the multiplier by discontinuous constant
or linear polynomials. The displacements are approximated on boundary edges by
linear functions.  Generalizations of the idea of weakly imposed symmetry to
other triangular elements, rectangular elements, and three space dimensions
were developed in \cite{Stenberg86}, \cite{Stenberg88},
\cite{Stenberg-mafelap} and \cite{Morley}.  In \cite{Stenberg88}, a family of
elements is developed in both two and three dimensions.  The lowest order
element in the family uses quadratics plus the curls of quartic bubble
functions in two dimensions or quintic bubble functions in three dimensions to
approximate the stresses, discontinuous linears to approximate the
displacements, and discontinuous quadratics to approximate the multiplier.  In
addition, a lower order method is introduced that approximates the stress by
piecewise linear functions augmented by the curls of cubic bubble functions
plus a cubic bubble times the gradient of local rigid motions.  The multiplier
is approximated by discontinuous piecewise linear functions and the
displacement by local rigid motions.  Morley \cite{Morley} extends PEERS to a
family of triangular elements, to rectangular elements, and to three
dimensions.  In addition, the multiplier is approximated by nonconforming
rather than continuous piecewise polynomials.

There is a close connection between mixed finite elements for linear
elasticity and discretization of an associated differential complex, the
elasticity complex, which will be introduced in \S~3 below. In fact, the
importance of this complex was already recognized in \cite{Arnold-Winther02},
where mixed methods for elasticity in two space dimensions were discussed.
The new ingredient here is that we utilize a constructive derivation of the
elasticity complex starting from the de Rham complex. This construction is
described in Eastwood \cite{Eastwood} and is based on the the
Bernstein--Gelfand--Gelfand resolution, cf.~\cite{BGG} and also \cite{Cap}.
By mimicking the construction in the discrete case, we are able to derive new
mixed finite elements for elasticity in a systematic manner from known
discretizations of the de Rham complex. As a result, we can construct new
elements both in two and three space dimensions which are significantly
simpler than those derived previously.  For example, we will construct stable
discretizations of the system \eqref{emixed-system} which only use piecewise
linear and piecewise constant functions, as illustrated in the figure below.
For simplicity, the entire discussion of the present paper will be given in
the three dimensional case.  A detailed discussion in two space dimensions can
be found in \cite{Arnold-Falk-Winther05_2}. Besides the methods discussed
here, we note that by slightly generalizing the approach of this paper, one
can also analyze some of the previously known methods mentioned above that are
also based on the weak symmetry formulation (see \cite{falk-cime1} for
details).

\begin{figure}[htb]
\centerline{\raise-.154in\hbox{\includegraphics[width=0.8in]{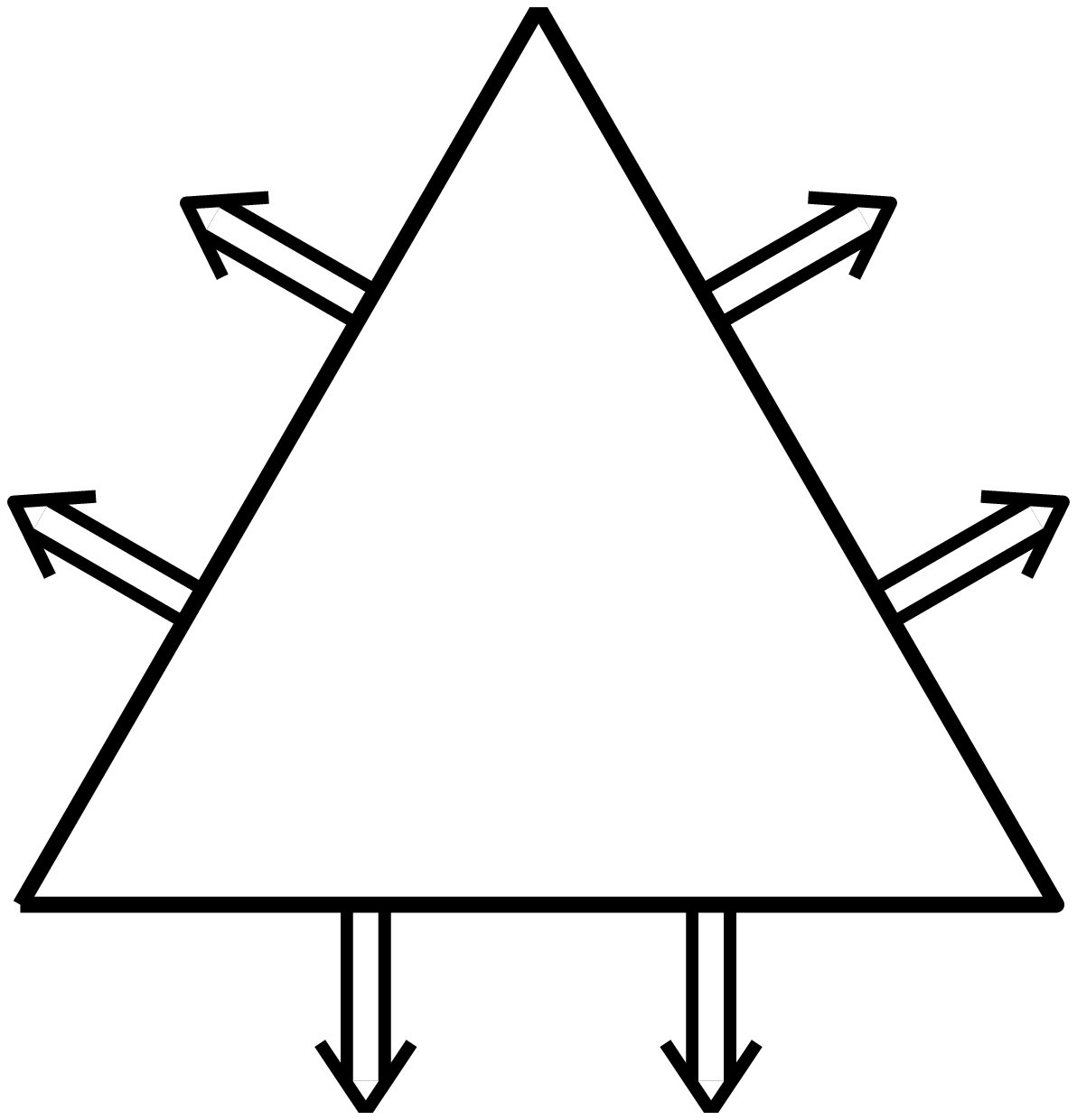}}
  \quad
  \includegraphics[width=0.8in]{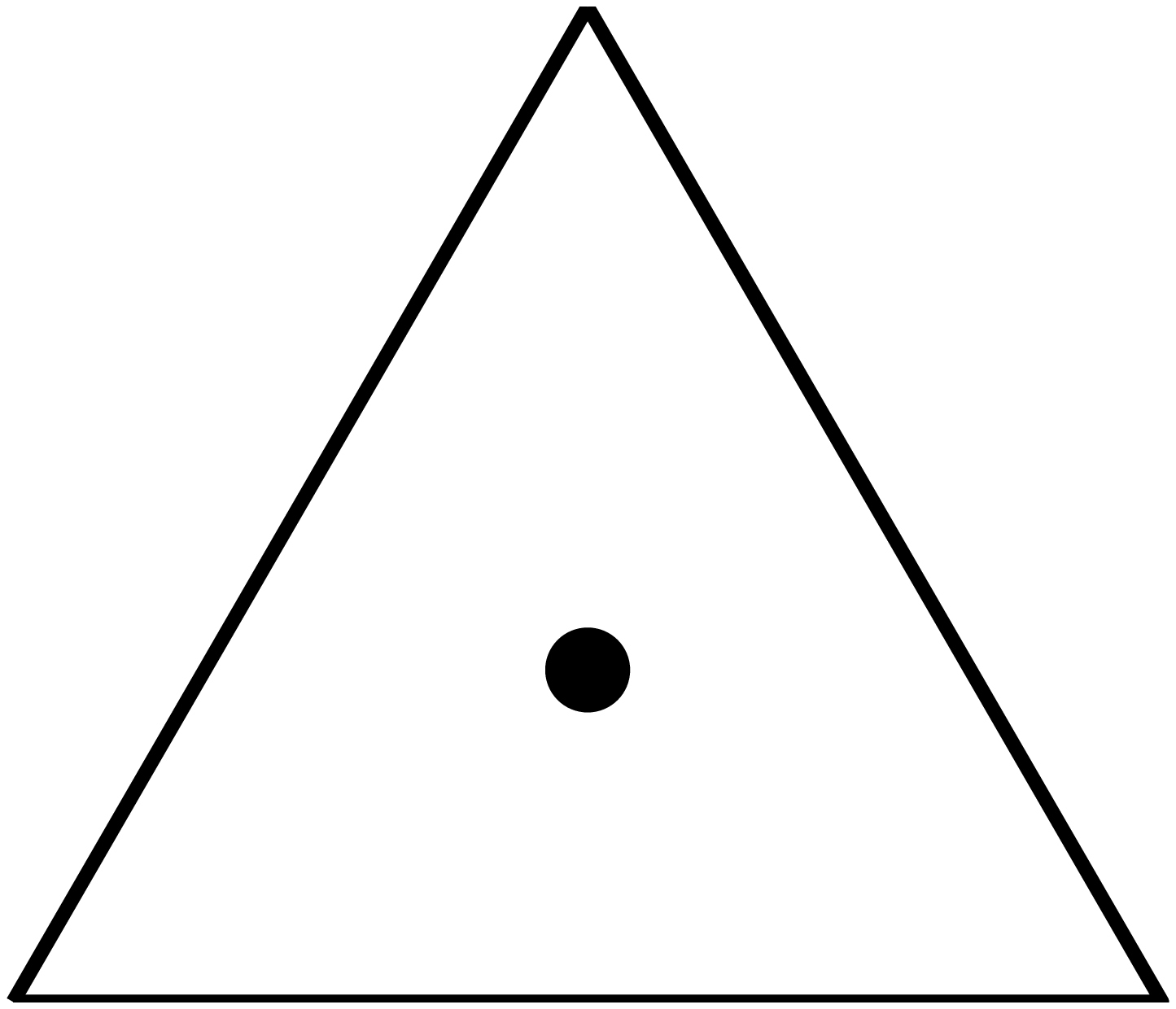}
  \quad
  \includegraphics[width=0.8in]{p0.eps}
  \qquad
  \raise-.256in\hbox{\includegraphics[width=0.8in]{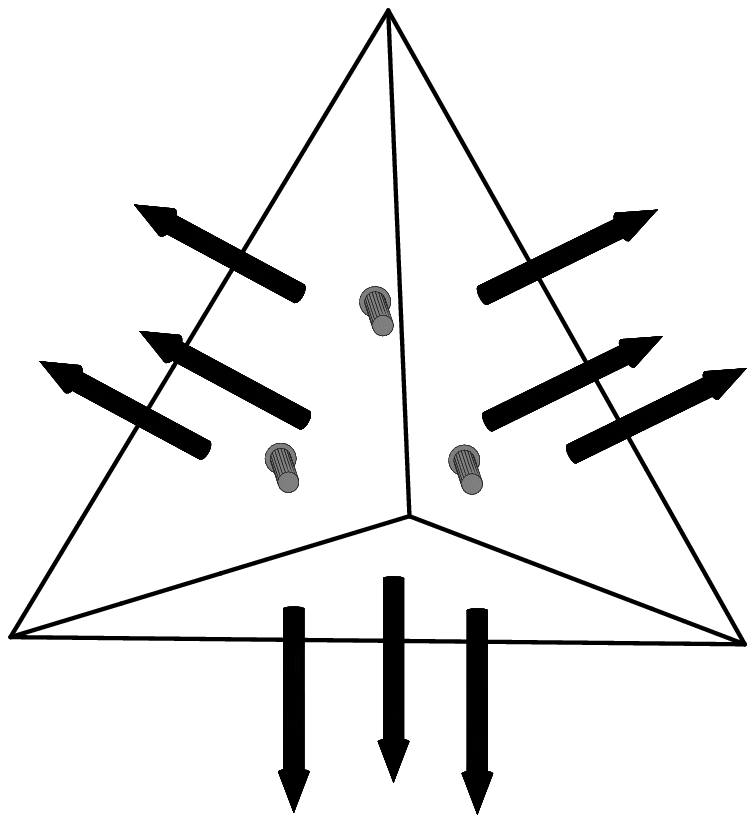}}
  \quad
  \includegraphics[width=0.8in]{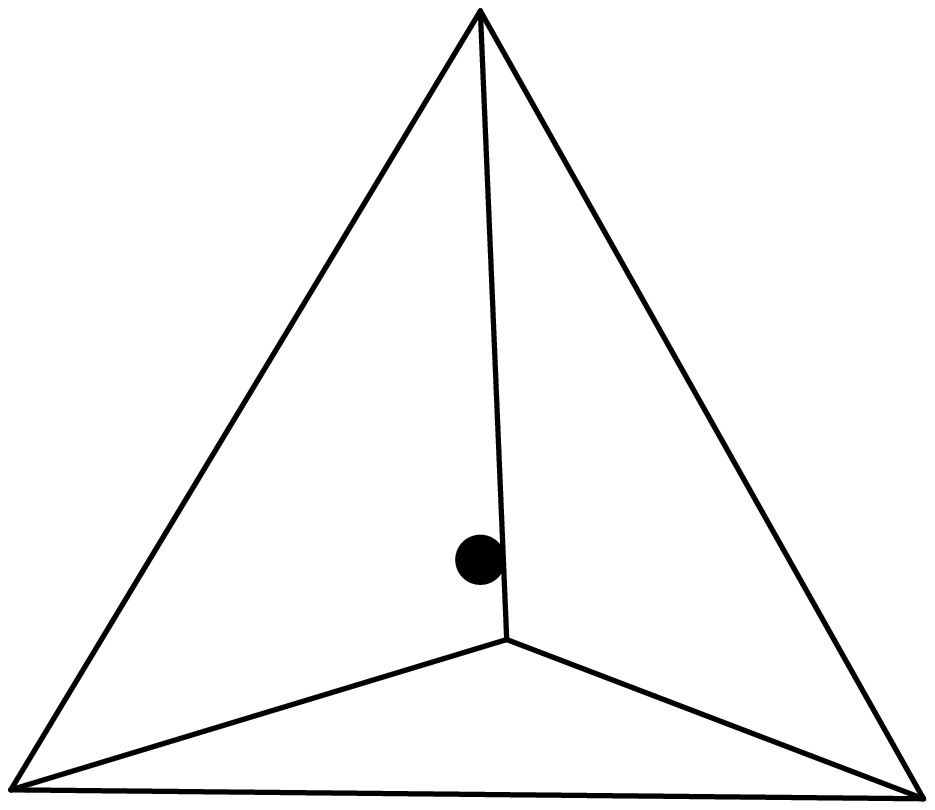}
  \quad
  \includegraphics[width=0.8in]{p03d.eps}}
\caption[]{Elements for the stress, displacement, and multiplier
in the lowest order case
in two dimensions and three dimensions.}
\end{figure}

An alternative approach to construct finite element methods for 
linear elasticity is to consider a pure displacement formulation.
Since the coefficient $A$ in \eqref{mixed-system} is invertible,
the stress $\sigma$ can be eliminated using the first equation in 
\eqref{mixed-system}, the stress--strain relation. This
leads to the second order
equation
\begin{equation}\label{second-order}
\div A^{-1} \eps u = f \quad \text{in } \Omega
\end{equation}
for the displacement $u$. A weak solution of this equation can be
characterized as the global minimizer of the energy functional
\[
\E(u)= \int_\Omega\bigl(\frac12
A^{-1}\eps u:\eps u + f\cdot u\bigr)\,dx
\]
over the Sobolev space $H^1_0(\Omega;\V)$. Here $H^1_0(\Omega;\V)$ denotes the
space of all square integrable vector fields on $\Omega$, with square
integrable derivatives, and which vanish on the boundary $\partial \Omega$. A
finite element approach based on this formulation, where we seek a minimum
over a finite element subspace of $H^1_0(\Omega;\V)$ is standard and discussed
in textbooks, (e.g., \cite{Ciarlet}). However, for more general models,
arising, for example, in viscoelasticity and plasticity
(cf.~\cite{Christensen}), the stress--strain relation is not local and an
elimination of the stress $\sigma$ is impossible.  For such models, a pure
displacement model is excluded, and a mixed approach seems to be an obvious
alternative. The construction of stable mixed elements for linear elasticity
is an important step in the construction of mixed methods for these more
complicated models. Another advantage of the mixed approach is that we
automatically obtain schemes which are uniformly stable in the incompressible
limit, i.e., as the Lam\'{e} parameter $\lambda$ tends to infinity. Since this
behavior of mixed methods is well known, we will not focus further on this
property here. A more detailed discussion in this direction can, for example,
be found in \cite{Arnold-Douglas-Gupta}.

An outline of the paper is as follows. In \S~2, we describe the
notation to be used, state our main result, and provide some
preliminary discussion on the relation between stability of mixed
finite element methods and discrete exact complexes.  In \S~3, we
present two complexes related to the two mixed formulations of
elasticity given by \eqref{mixed-system_w} and \eqref{emixed-system}. 
In \S~4, we introduce the framework of differential forms and show how
the elasticity complex can be derived from the de Rham complex. In
\S~5, we derive discrete analogues of the elasticity complex beginning
from discrete analogues of the de Rham complex and identify the
required properties of the discrete spaces necessary for this
construction.  This procedure is our basic design principle.  In \S~6,
we apply the construction of the preceding section to specific discrete
analogues of the de~Rham complex to obtain a family of discrete
elasticity complexes.  In \S~7 we use this family to construct stable
finite element schemes for the approximation of the mixed formulation
of the equations of elasticity with weakly imposed symmetry.  Finally,
in \S~8, we show how a slightly more complicated procedure leads to a
simplified elasticity element.

\section{Notation, statement of main results, and preliminaries}
\label{sec:prelim}
We begin with some basic notation and hypotheses.  We continue to denote by
$\V=\R^3$ the space of $3$-vectors, by
$\M$ the space of $3\times 3$ real matrices, and by
$\Sym$ and $\K$ the subspaces of symmetric and skew symmetric matrices,
respectively.  The operators $\sym:\M\to\Sym$ and $\skw:\M\to\K$ denote
the symmetric and skew symmetric parts, respectively.
Note that an element of the space $\K$ can be identified with its
\emph{axial vector} in $\V$ given by the map
$\vect : \K \to \V$:
\begin{equation*}
\vect \begin{pmatrix} 0 & -v_3 & v_2 \\
v_3 & 0 & -v_1 \\ -v_2 & v_1 & 0 \end{pmatrix}
= \begin{pmatrix} v_1 \\ v_2 \\ v_3 \end{pmatrix},
\end{equation*}
i.e., $\vect^{-1}(v)w= v\x w$ for any vectors $v$ and $w$.

We assume that $\Omega$ is a  domain in $\R^3$ with
boundary $\partial \Omega$.  We shall use 
the standard function spaces, like the
Lebesgue space $L^2(\Omega)$ and the Sobolev space $H^s(\Omega)$.  For
vector-valued functions, we include the range space in the notation
following a semicolon, so $L^2(\Omega;\X)$ denotes the space of square
integrable functions mapping $\Omega$ into a normed vector space $\X$.  The
space $H(\div,\Omega;\V)$ denotes the subspace of (vector-valued)
functions in $L^2(\Omega;\V)$ whose divergence belongs to
$L^2(\Omega)$.  Similarly, $H(\div,\Omega;\M)$ denotes the subspace of
(matrix-valued) functions in $L^2(\Omega;\M)$ whose
divergence (by rows) belongs to $L^2(\Omega;\V)$.

Assuming that $\X$ is an inner product space, then $L^2(\Omega;\X)$
has a natural norm and inner product, which will be denoted by
$\|\,\cdot\,\|$ and $(\,\cdot\,,\,\cdot\,)$, respectively.
For a Sobolev space $H^s(\Omega;\X)$, we denote the norm by
$\|\,\cdot\,\|_s$ and for $H(\div, \Omega;\X)$, the norm is denoted
by $\|v \|_{\div} := (\|v\|^2 + \|\div v\|^2)^{1/2}$.
The space
$\P_k(\Omega)$ denotes the space of
polynomial functions on $\Omega$ of total degree $\le k$.  Usually
we abbreviate this to just $\P_k$.

In this paper we shall consider mixed finite element approximations
derived from \eqref{emixed-system}.  These schemes  take the form:

Find $(\sigma_h, u_h, p_h) \in \Sigma_h \times V_h \times Q_h$ such that
\begin{equation}\label{emixed-system-h}
\begin{array}{lrll}
\int_\Omega(A\sigma_h:\tau + \div\tau\cdot u_h + \tau : p_h)\,dx &=& 0,
\quad &\tau \in \Sigma_h,\\[0.2cm]
\int_\Omega \div\sigma_h \cdot v \,dx  &=& 
\int_\Omega f \cdot v\,dx, \quad
&v \in V_h,\\[0.2cm]
\int_\Omega \sigma_h :q \,dx  &=& 0, \quad
&q \in Q_h,
\end{array}
\end{equation}
where now  $\Sigma_h\subset \Hdivn$, $V_h\subset
L^2(\Omega;\V)$ and $Q_h \in L^2(\Omega; \K)$.

Following the general theory of
mixed finite element methods, cf.~\cite{Brezzi,Brezzi-Fortin}, 
the stability of the
saddle--point system \eqref{emixed-system-h} is ensured by the following
conditions:
\begin{itemize}
\setlength{\itemindent}{20pt}
\item[(A$1$)] 
$\norm{\tau}^2_{\div} \le c_1(A\tau, \tau)$ whenever
$\tau\in\Sigma_h$ satisfies
$(\div\tau, v) =0\quad \forall v\in V_h$, 
and $(\tau, q) =0 \quad \forall q\in Q_h$, 
\smallskip
\item[\quad(A$2$)] for all nonzero $(v,q) \in V_h \x Q_h$,
there exists nonzero $\tau\in\Sigma_h$ with $(\div\tau, v) + (\tau, q) \ge
c_2 \norm{\tau}_{\div}(\norm{v} + \norm{q})$,
\end{itemize}
where $c_1$ and $c_2$ are positive constants independent of $h$.

The main result of this paper, given in Theorem~\ref{thm1-A2}, is to 
construct a new family of stable
finite element spaces $\Sigma_h$, $V_h$, $Q_h$ that satisfy the
stability conditions (A$1$) and (A$2$).  We shall show that for $r \ge
0$, the choices of the N\'ed\'elec second family of  H(div) elements of
degree $r+1$ for $\Sigma_h$, cf.~\cite{Nedelec2}, and of discontinuous
piecewise polynomials of degree $r$ for $V_h$ and $Q_h$ provide a
stable finite element approximation.  In contrast to the previous work
described in the introduction, no stabilizing bubble functions are
needed; nor is interelement continuity imposed on the multiplier.
In \S~8 we also discuss a somewhat 
simpler lowest order element ($r=0$) in
which the local stress space is a strict subspace of the full space of linear
matrix fields.

Our approach to the construction of stable mixed elements for
elasticity is motivated by the success in developing stable mixed
elements for steady heat conduction (i.e., the Poisson problem) based on
discretizations of the de~Rham complex. We recall (see, e.g.,
\cite{Arnold-Falk-Winther05_1}) that there is a close connection
between the construction of such elements and discretizations of
the de~Rham complex
\begin{equation}\label{deRham1}
\begin{CD}
\R  @.\,\hookrightarrow\,@. C^\infty(\Omega) @>\grad>> 
C^\infty(\Omega;\V)@>\curl>>
C^\infty(\Omega; \V) @>\div>> C^\infty(\Omega) @>>> 0.
\end{CD}
\end{equation}
More specifically, a key to the construction and analysis of stable
mixed elements is a commuting diagram of the form
\begin{equation}\label{diagram4}
\begin{CD}
\R  @.\,\hookrightarrow\,@. C^\infty(\Omega) @>\grad>> 
C^\infty(\Omega;\V)@>\curl>>
C^\infty(\Omega; \V) @>\div>> C^\infty(\Omega) @>>> 0\\
@. @. @VV\Pi^1_hV @VV\Pi_h^cV @VV\Pi^d_hV  @VV\Pi^0_hV\\
\R  @.\,\hookrightarrow\,@. W_h @>\grad>> U_h @>\curl>>
V_h @>\div>> Q_h @>>> 0.
\end{CD}
\end{equation}
Here, the spaces
$V_h\subset H(\div)$ and $Q_h\subset L^2$ are the finite element spaces
used to discretize the flux and temperature fields, respectively.  The 
spaces $U_h\subset H(\curl)$ and $W_h\subset H^1$ are additional finite
element spaces, which can be found for all well-known stable element
choices.  The bottom row of the diagram is a discrete de~Rham complex,
which is exact when the de~Rham complex is (i.e., when the domain is
contractible).  The vertical operators are projections determined by
the natural degrees of freedom of the finite element spaces.
As pointed out in \cite{Arnold-Falk-Winther05_1}, there are many such
discretizations of the de~Rham complex.

A diagram analogous to \eqref{diagram4}, but with the de~Rham complex
replaced by the \emph{elasticity complex} defined just below, will be
crucial to our construction of stable mixed elements for elasticity.
Discretization of the elasticity complex also gives insight into the
difficulties of constructing finite element approximations of the mixed
formulation of elasticity with strongly imposed symmetry,
cf.~\cite{Arnold-Falk-Winther05_2}.

\section{The elasticity complex}\label{sec:elas-complex}
We now proceed to a description of two elasticity complexes,
corresponding to strongly or weakly imposed symmetry of the stress tensor.
In the case of strongly imposed symmetry, relevant to the
mixed elasticity system \eqref{mixed-system_w}, 
the characterization of the divergence-free symmetric matrix fields
will be needed. In order to give such a characterization, define
$\curl : C^\infty(\Omega; \M) \to C^\infty(\Omega; \M)$
to be the differential operator defined by taking $\curl$ of each row of
the matrix. Then define a second order differential operator 
$J : C^\infty(\Omega;\Sym) \to C^\infty(\Omega;\Sym)$ by
\begin{equation}\label{curlcurl-op1}
J \tau = \curl (\curl \tau)^T, \quad \tau\in C^\infty(\Omega;\Sym).
\end{equation}
It is easy to check that $\div \circ J =0$ and that $J \circ \eps= 0$.
In other words, 
\begin{equation}\label{elas-seq1}
\begin{CD}
\T @.\,\hookrightarrow\,@. C^\infty(\V) @>\eps>>
C^\infty(\Sym)
@>J>>
C^\infty(\Sym) @>\div>> C^\infty(\V) @>>> 0,
\end{CD}
\end{equation}
is a complex. 
Here the dependence of the domain $\Omega$ is suppressed, i.e.,
$C^\infty(\Sym)= C^\infty(\Omega; \Sym)$, and
$\T = \T(\Omega)$ denotes the six dimensional space of infinitesimal rigid motions
on $\Omega$, i.e., functions of the form $x\mapsto a+Bx$ with $a\in\V$ and $B\in\K$.
In fact, when  $\Omega$ is contractible, then
\eqref{elas-seq1}
is an exact sequence, a fact which will follow from the discussion below.
The complex \eqref{elas-seq1}
will be referred to as the elasticity complex.

A natural approach to the construction of stable
mixed finite elements for elasticity would be to extend  the
complex \eqref{elas-seq1} to a complete commuting diagram of the form
\eqref{diagram4}, where \eqref{elas-seq1} is the top row and the bottom
row is a discrete analogue. However, due to the pointwise symmetry
requirement on the discrete stresses, this construction requires
piecewise polynomials of high order. For the corresponding problem in
two space dimensions, such a complex was proposed in
\cite{Arnold-Winther02} with a piecewise cubic stress space, cf.~also
\cite{Arnold-Falk-Winther05_2}. An analogous complex was derived in the
three dimensional case in \cite{Arnold-Awanou-Winther}.  It uses a
piecewise quartic space, with 162 degrees of freedom on each
tetrahedron for the stresses.

We consider the formulation based on weakly imposed symmetry of the
stress tensor, i.e., the mixed system \eqref{emixed-system}.  Then
the relevant complex is, instead of \eqref{elas-seq1}, 
\begin{equation}\label{elas-seq3}
\begin{CD}
\T' @.\,\hookrightarrow\,@. C^\infty(\V\x\K) @>(\grad, -I)>>
C^\infty(\M)
@>J>>
C^\infty(\M) @>(\div, \skw)^T >> C^\infty(\V\x\K) @>>> 0.
\end{CD}
\end{equation}
Here,
\[
\T' = \{\,(v,\grad v)\,|\,v\in\T\,\},
\]
and $J : C^\infty(\Omega;\M) \to C^\infty(\Omega;\M)$ denotes the extension of 
the operator defined on $C^\infty(\Omega;\Sym)$
by \eqref{curlcurl-op1} such that $J \equiv 0$ on 
$C^\infty(\Omega;\K)$.  We remark that $J$ may be written
\begin{equation}\label{curlcurl-op2}
J\tau = \curl \Xi^{-1}\curl \tau,
\end{equation}
where $\Xi :\M \to \M$ is the algebraic operator
\begin{equation}\label{S-oper1}
\Xi\mu = \mu^T - \tr(\mu)\delta, \qquad \Xi^{-1}\mu = \mu^T - 
\frac{1}{2}\tr(\mu)\delta,
\end{equation}
with $\delta$ the identity matrix.
Indeed, if $\tau$ is symmetric, then $\curl \tau$ is trace
free, 
and therefore the definition \eqref{curlcurl-op2} reduces to 
\eqref{curlcurl-op1} on $C^\infty(\Omega;\Sym)$.  On the other
hand, if $\tau$ is skew with axial vector $u$, then
$\curl\tau= -\Xi\grad u$, and so $\curl \Xi^{-1}\curl\tau=0$.

Observe that there is a close connection between
\eqref{elas-seq1} and \eqref{elas-seq3}. In fact, \eqref{elas-seq1}
can be derived from \eqref{elas-seq3} by performing a projection step.
To see this, consider the diagram
\begin{equation}\label{symmetrization}
\begin{CD}
\T'@.\,\hookrightarrow\,@. C^\infty(\V\x\K) @>(\grad, -I)>>
C^\infty(\M)
@>J>>
C^\infty(\M) @>(\div, \skw)^T >> C^\infty(\V\x\K) @>>> 0
\\
 @. @. @VV\pi_0V @VV\pi_1V
 @VV\pi_2V  @VV\pi_3V
\\
\T @.\,\hookrightarrow\,@. C^\infty(\V) @>\eps>>
C^\infty(\Sym)
@>J>>
C^\infty(\Sym)  @>\div>> C^\infty(\V)  @>>> 0,
\end{CD}
\end{equation}
where the projection operators $\pi_k$ are defined by
\[
\pi_0(u,q) = u, \quad \pi_1(\sigma) = \pi_2(\sigma) = \sym(\sigma),
\quad \pi_3(u,q) = u - \div q.
\]
We may identify $C^\infty(\V)$ with a subspace of
$C^\infty(\V\x\K)$, namely,
\[
\{(u,q): \, u \in C^\infty(\V), q = \skw(\grad u)\}.
\]
Under this identification, $\T\subset C^\infty(\V)$
corresponds to $\T'\subset C^\infty(\V\x\K)$. 
We identify the $C^\infty(\V)$ on the
right with a different subspace of $C^\infty(\V\x\K)$, namely,
\[
\{(u,q): \, u \in C^\infty(\V), q = 0\}.
\]
With these identifications, the bottom row is a subcomplex of the top row,
and the operators $\pi_k$ are all projections. Furthermore, the diagram
commutes.  It follows easily that the exactness of the upper row
implies exactness of the bottom row.

In the next section, we shall discuss these complexes further. In particular,
we show the elasticity complex with weakly imposed symmetry, i.e.,
\eqref{elas-seq3}, follows from the de Rham complex \eqref{deRham1}. Hence,
as a consequence of the discussion above, both \eqref{elas-seq1} and
\eqref{elas-seq3} will follow from \eqref{deRham1}.

\section{From the de~Rham to the elasticity complex}\label{sec:BGG}
The main purpose of this section is to demonstrate the connection
between the de Rham complex \eqref{deRham1} and the elasticity
complex \eqref{elas-seq3}.  In particular, we show that whenever
\eqref{deRham1} is exact, \eqref{elas-seq3} is exact.
This section serves as an introduction to a corresponding
construction of a discrete elasticity complex, to be given in the next
section.  In the following section,
the discrete complex will be used to construct
stable finite elements for the system \eqref{emixed-system}.

The de Rham complex  \eqref{deRham1} is most clearly stated in terms of
differential forms.  Here we briefly recall the definitions and
properties we will need.  We use a completely coordinate-free approach.
For a slightly more expanded discussion and the expressions in coordinates
see, e.g., \cite[\S~4]{Arnold-Falk-Winther05_1}.
We let $\Lambda^k$ denote the space of smooth differential
$k$-forms on $\Omega$, i.e.~$\Lambda^k = \Lambda^k(\Omega)
= C^\infty(\Omega; \text{Alt}^k\V)$,
where $\text{Alt}^k\V$ denotes the vector space of alternating $k$-linear maps
on $\V$. If $\omega \in \Lambda^k$ we let $\omega_x \in
\text{Alt}^k\V$
denote $\omega$ evaluated at $x$, i.e., we use subscripts to indicate
the spatial dependence.

Using the inner product on $\text{Alt}^k\V$ inherited from the
inner product on $\V$ (see equation (4.1) of \cite[\S~4]{Arnold-Falk-Winther05_1}),
we may also define the Hilbert space 
$L^2\Lambda^k(\Omega)=L^2(\Omega;\text{Alt}^k\V)$ of square integrable differential forms
with norm denoted by $\|\,\cdot\,\|$, and also the $m$th order Sobolev space
$H^m\Lambda^k(\Omega)=H^m(\Omega;\text{Alt}^k\V)$,
consisting of square integrable $k$-forms for which the norm
\[
\|\omega\|_m := \bigl( \sum_{|\alpha|\le m} 
\|\partial^\alpha \omega\|^2\bigr)^{1/2}
\]
is finite (where the sum is over multi-indices of degree at most $m$).

%We recall the usual wedge product
%$\wedge:\Lambda^k\x\Lambda^l\to\Lambda^{k+l}$ given by
%\begin{equation*}
%(\omega\wedge\mu)(v_1,\ldots,v_{k+l})=\sum
%(\operatorname{sign}\sigma)\omega(v_{\sigma_1},\ldots,v_{\sigma_k})
%\mu(v_{\sigma_{k+1}},\ldots,v_{\sigma_{k+l}}),\quad 
%\omega\in\Lambda^k, \mu\in\Lambda^l, v_i\in V,
%\end{equation*}
%where the sum is over the set of all permutations of 
%$\{1,\ldots,k+l\}$, for which $\sigma_1<\sigma_2<\cdots<\sigma_k$
%and $\sigma_{k+1}<\sigma_{j+2}<\cdots<\sigma_{k+l}$.  We write
%$e_i$ and $dx^i$, $i=1,2,3$, for the usual basis and dual
%basis elements of $\V$.
%Then each $\omega \in \Lambda^k$ can be represented as 
%\begin{equation}\label{k-form}
%\omega_x = \sum_{1 \le i_1<i_2<\ldots <i_k \le 3} f_{i_1 \ldots i_k}(x) 
%dx^{i_1}\wedge \ldots \wedge dx^{i_k} = \sum_{I} f_I(x) dx^I
%\end{equation}
%with coefficients $f_I \in C^\infty(\Omega;\R)$ for each $k$-tuple $I$.
%Note that,
%\[
%dx^{i_1}\wedge \ldots \wedge dx^{i_k}(v_1,\ldots,v_k)
%= \text{det}(v_{{i_j},m})_{j,m=1}^k, \quad 1 \le i_1<i_2<\ldots <i_k \le 3,
%\]
%where $v_m = (v_{1,m}, v_{2,m}, v_{3,m})^T$.
Thus, $0$-forms are scalar functions and $1$-forms
are covector fields.   We will not emphasize the distinction between
vectors and covectors, since, given the inner product in $\V$, we may identify
a $1$-form $\omega$ with the vector field $v$ for which
$\omega(p)=v\cdot p$, $p\in\V$.  In the $3$-dimensional case,
we can identify a $2$-form $\omega$ with a vector field $v$ and a $3$-form
$\mu$ with a scalar field $c$ by
\[
\omega(p,q)=v\cdot p\x q, \quad \mu(p,q,r) = c(p \x q \cdot r),
\quad p,q,r\in\V.
\]

The exterior derivative $d=d_k: \Lambda^k \to\Lambda^{k+1}$
is defined by
\begin{equation}\label{extd}
d\omega_x(v_1,\ldots,v_{k+1})=\sum_{j=1}^{k+1} (-1)^{j+1}
\partial_{v_j}\omega_x(v_1,\ldots,\hat v_j,\ldots,v_{k+1}),\quad
\omega\in\Lambda^k, v_1,\ldots,v_{k+1}\in\V,
\end{equation}
where the hat is used to indicate a suppressed argument
and $\partial_{v}$ denotes the directional derivative in the direction
of the vector $v$.  It is useful to define
\[
H\Lambda^k = \{\, \omega\in L^2(\Omega;\text{Alt}^k\V)\,|\, d\omega
\in L^2(\Omega;\text{Alt}^{k+1}\V)\,\},
\]
with norm given by $\|\omega\|_{H\Lambda}^2=\|\omega\|^2+\|d\omega\|^2$.
%In traditional notation, the exterior derivative
%is
%\[
%d \sum_{i_1<\cdots<i_k} f_{i_1\ldots i_k} dx^{i_1}\wedge\cdots\wedge dx^{i_k}
%=\sum_j \sum_{i_1<\cdots<i_k} \frac{\partial f_{i_1\ldots i_k}}{\partial x^j}
%dx^j\wedge
%dx^{i_1}\wedge\cdots\wedge dx^{i_k}.
%\]
Using the identifications given above, the $d_k$ correspond to 
$\grad$, $\curl$, and $\div$ for $k=0,1,2$, respectively, and
the $H\Lambda^k$ correspond to $H^1$, $H(\curl)$, $H(\div)$, and,
for $k=3$, $L^2$.
%We also recall the Leibniz rule
%\[
%d(\omega \wedge \mu) = d\omega \wedge \mu + (-1)^k \omega \wedge d\mu,
%\quad \omega \in \Lambda^k, \mu \in \Lambda^\ell.
%\]

The de Rham complex \eqref{deRham1} can then be written
\begin{equation}\label{deRham-form1}
\begin{CD}
\R  @.\,\hookrightarrow\,@. \Lambda^0 @>d>> \Lambda^1 @>d>>
\Lambda^2@>d>> \Lambda^3@>>> 0.
\end{CD}
\end{equation}
It is a complex since $d\circ d=0$.

A differential $k$-form $\omega$ on $\Omega$ may be restricted to a differential
$k$-form on any submanifold
$M\subset\bar\Omega$: at each point of $M$ the restriction of $\omega$ is
an alternating linear form on tangent vectors.  Moreover if $\dim M=k$,
the integral $\int_M \omega$ is defined.

%Next, we define for $\omega \in \Lambda^k(\Omega)$ and $v_1, \cdots, v_k$
%tangent vectors to $\Gamma$, the trace $\Tr \omega(v_1, \cdots, v_k)$ $=
%\omega(v_1, \cdots, v_k)$ for all $x \in \Gamma$ and the boundary exterior
%derivative $d_{\Gamma}: \Lambda_{\Gamma}^k \to \Lambda_{\Gamma}^{k+1}$ 
%satisfying $d_{\Gamma} \Tr \omega = \Tr d \omega$. 
%We then have the commuting diagram
%\begin{equation}\label{deRham-boundary}
%\begin{CD}
% \R  @.\,\hookrightarrow\,@. \Lambda^0 @>d>> \Lambda^1 @>d>>
%\Lambda^2@>d>> \Lambda^3@>>> 0
%\\
%@VV id V @. @VV \Tr V @VV \Tr V @VV \Tr V 
%\\
%\R  @.\,\hookrightarrow\,@. \Lambda_{\Gamma}^0 @>d_{\Gamma}>>
% \Lambda_{\Gamma}^1 @>d_{\Gamma}>>
%\Lambda_{\Gamma}^2@> >> 0
%\end{CD}.
%\end{equation}
%If $\mu_x = g(x) dx^1\wedge dx^2\wedge dx^3$ is a $3$-form on
%$\Omega$
%then $\int_\Omega \mu$ is simply the integral of the function $g$.
%Similarly, if $\omega$ is a $2$-form, identified with the vector field 
%$f$, then
%\[
%\int_{\partial \Omega} \Tr \omega = \int_{\partial \Omega} f\cdot n\,
%d\mathfrak H_2,
%\]
%where $n$ is the outward unit normal vector and $\mathfrak H_2$
%is the $2$-dimensional Hausdorff measure.
%Hence, Stokes theorem takes the form
%\[
%\int_\Omega d\omega = \int_{\partial \Omega} \Tr \omega.
%\]

If $\X$ is a vector space, then $\Lambda^k(\X) = \Lambda^k(\Omega;\X)$
refers to the $k$-forms with values in $\X$,
i.e., $\Lambda^k(\X) = C^\infty(\Omega; \text{Alt}^k(\V;\X))$, where
$\text{Alt}^k(\V;\X)$ are alternating $k$-linear forms on $\V$
with values in $\X$.  Given an inner product on $\X$, we
obtain an inner product on $\Lambda^k(\X)$.
Obviously the corresponding complex
\begin{equation}\label{deRham-form2}
\begin{CD}
\X  @.\,\hookrightarrow\,@. \Lambda^0(\X) @>d>> \Lambda^1 (\X)@>d>>
\Lambda^2(\X) @>d>> \Lambda^3(\X)@>>> 0,
\end{CD}
\end{equation}
is exact whenever the de~Rham complex is.

We now construct the elasticity complex as a subcomplex of a complex
isomorphic to the de~Rham complex with values in the six-dimensional
vector space $\W:=\K\x\V$.
First, for any $x\in\R^3$ we define $K_x:\V\to\K$ by $K_x v =2\skw(xv^T)$.
We then define an operator $K:\Lambda^k(\Omega;\V)\to\Lambda^k(\Omega;\K)$
by
\begin{equation}\label{defk}
(K\omega)_x(v_1,\ldots,v_k) = K_x[\omega_x(v_1,\ldots,v_k)].
\end{equation}
Next, we define an isomorphism
$\Phi:\Lambda^k(\W)\to \Lambda^k(\W)$ by
\[
\Phi(\omega,\mu) = (\omega + K\mu, \mu),
\]
with inverse given by 
\[
\Phi^{-1}(\omega,\mu) = (\omega - K\mu, \mu).
\]
Next, define the operator $\A : \Lambda^k(\W) \to \Lambda^{k+1}(\W)$
by $\A = \Phi d \Phi^{-1}$. Inserting the isomorphisms $\Phi$ in the
$\W$-valued de~Rham sequence, we obtain a complex
\begin{equation}\label{deRham-form3}
\begin{CD}
\Phi(\W)  @.\,\hookrightarrow\,@. \Lambda^0(\W) @>\A>> \Lambda^1 (\W)@>\A>>
\Lambda^2(\W) @>\A>> \Lambda^3(\W)@>>> 0,
\end{CD}
\end{equation}
which is exact whenever the de~Rham complex is.

The operator $\A$ has a simple form. Using the definition of 
$\Phi$, we obtain for $(\omega,\mu) \in \Lambda^k(\W)$,
\begin{equation*}
\A(\omega,\mu) = \Phi \circ d (\omega - K\mu,\mu) = 
\Phi(d\omega - d K\mu,d\mu)
= (d\omega - S \mu, d\mu),
\end{equation*}
where $S= S_k:\Lambda^k(\V) \to \Lambda^{k+1}(\K)$, $k=0,1,2$, is given by
$S=dK - Kd$.  Using
the definition \eqref{extd} of the exterior derivative, the
definition \eqref{defk} of $K$, and the Leibniz rule 
\begin{equation}\label{leibniz}
d(\omega \wedge \mu) = d\omega \wedge \mu + (-1)^k \omega \wedge d\mu,
\quad \omega \in \Lambda^k,~\mu \in \Lambda^\ell,
\end{equation}
we obtain
\begin{equation}\label{defS}
(S\omega)(v_1,\ldots,v_{k+1}) = \sum_{j=1}^{k+1}(-1)^{j+1}
K_{v_j}[\omega(v_1,\ldots,\hat v_j,\ldots v_{k+1})],\quad
\omega\in \Lambda^k(\Omega;\V).
\end{equation}
Note that the operator $S$ is purely algebraic, and independent of $x$.

Since $d^2 =  0$ we have 
\[
dS = d^2 K - dKd = - (dK - Kd)d
\]
or
\begin{equation}\label{dS=-Sd}
dS = - Sd.
\end{equation}
Noting that
\begin{multline*}
(S_1\mu)(v_1,v_2)=K_{v_1}[\mu(v_2)]-K_{v_2}[\mu(v_1)]=2\skw
[v_1\mu(v_2)^T-v_2\mu(v_1)^T],
\\
\mu\in\Lambda^1(\Omega;\V),~v_1,v_2\in\V,
\end{multline*}
we find, using the identity 
\begin{equation}\label{cross-skew}
a\x b=-2\vect\skw ab^T,
\end{equation}
that $S_1$ is invertible with
\begin{multline*}
(S_1^{-1}\omega)(v_1)\x v_2\cdot v_3=\frac12[
  \vect\bigl(\omega(v_2,v_3)\bigr)\cdot v_1
- \vect\bigl(\omega(v_1,v_2)\bigr)\cdot v_3
+ \vect\bigl(\omega(v_1,v_3)\bigr)\cdot v_2
],
\\
\omega\in\Lambda^2(\Omega;\K),v_1, v_2, v_3\in\V.
\end{multline*}

We now define the desired subcomplex.  Define
\[
\Gamma^1 = \{\,(\omega, \mu) \in \Lambda^1(\Omega;\W)\,|\, d\omega = S_1\mu \,\},
\quad
\Gamma^2 = \{\,(\omega, \mu) \in \Lambda^2(\Omega;\W)\,|\, \omega = 0\,\},
\]
with projections $\pi^1:\Lambda^1(\Omega;\W)\to\Gamma^1$ and
$\pi^2:\Lambda^2(\Omega;\W)\to\Gamma^2$ given by
\[
\pi^1(\omega, \mu) = (\omega, S_1^{-1}d \omega), \quad
\pi^2(\omega, \mu) = (0 , \mu + d S_1^{-1} \omega).
\]
Using \eqref{dS=-Sd}, it is straightforward to check that
$\A$ maps $\Lambda^0(\W)$ into $\Gamma^1$ and $\Gamma^1$ into
$\Gamma^2$, and that the diagram
\begin{equation}\label{projection1}
\begin{CD}
\Phi(\W)  @.\,\hookrightarrow\,@. \Lambda^0(\W) @>\A>> \Lambda^1 (\W)@>\A>>
\Lambda^2(\W) @>\A>> \Lambda^3(\W)@>>> 0\\
@. @. @VVidV @VV\pi^1V @VV\pi^2V  @VVidV\\
\Phi(\W) @.\,\hookrightarrow\,@. \Lambda^0(\W) @>\A>> \Gamma^1@>\A>>
\Gamma^2@>\A>> \Lambda^3(\W)@>>> 0
\end{CD}
\end{equation}
commutes, and therefore the subcomplex in the bottom row is exact when
the de~Rham complex is.  This
subcomplex is, essentially, the elasticity complex.
Indeed, by identifying elements $(\omega,\mu) \in \Gamma^1$
with $\omega \in \Lambda^1(\K)$, and elements $(0,\mu) \in \Gamma^2$
with $\mu \in \Lambda^2(\V)$, the subcomplex becomes
\begin{equation}\label{elas-seq4}
\begin{CD}
\Phi(\W)  @.\,\hookrightarrow\,@. \Lambda^0(\K\x\V) @>(d_0,-S_0)>> \Lambda^1 (\K)@>d_1\circ
S_1^{-1}\circ d_1>>
\Lambda^2(\V) @>(-S_2,d_2)^T>> \Lambda^3(\K\x\V)@>>> 0.
\end{CD}
\end{equation}
This complex may be identified with 
\eqref{elas-seq3}. As an initial step of this identification we
observe that the algebraic operator $\Xi : C^\infty(\M)\to
C^\infty(\M)$
appearing in \eqref{elas-seq3} via \eqref{curlcurl-op2} and the operator $S_1 : \Lambda^1(\V)
\to \Lambda^2(\K)$ are connected by the identity
\begin{equation}\label{identification}
\Xi = \Upsilon_2^{-1}S_1\Upsilon_1,
\end{equation}
where $\Upsilon_1: C^\infty(\M) \to \Lambda^1(\V)$ and 
$\Upsilon_2: C^\infty(\M) \to \Lambda^2(\K)$ are given by
$\Upsilon_1F(v) =Fv$ and $\Upsilon_2 F(v_1,v_2) = \vect^{-1}F(v_1 \x
v_2)$ for $F \in C^\infty(\M)$.
In fact, using \eqref{cross-skew}, we have for any $v_1, v_2 \in \V$
\begin{align*}
S_1 \Upsilon_1F(v_1,v_2) &= 2\skw[v_1(Fv_2)^T -v_2(Fv_1)^T]\\
&= \vect^{-1}(v_2 \x Fv_1 - v_1 \x Fv_2).
\end{align*}
On the other hand,
\[
\Upsilon_2 \Xi F(v_1,v_2) = \vect^{-1}[\Xi F(v_1 \x v_2)],
\]
and hence \eqref{identification} follows from the algebraic identity
\[
\Xi F (v_1 \x v_2) = v_2 \x Fv_1 - v_1 \x Fv_2,
\]
which holds for any $F \in \M$.

We may further identify the four spaces of fields in
\eqref{elas-seq3} with the corresponding spaces of forms in 
\eqref{elas-seq4} in a natural way:
\begin{itemize}
\item $(u,p)\in C^\infty(\V\x\K)\sim
(\vect^{-1}u,\vect p)\in\Lambda^0(\K\x\V)$.
\item $F\in C^\infty(\M) \sim \omega\in\Lambda^1(\K)$
given by $\omega(v)=\vect^{-1}(Fv)$.
\item $F\in C^\infty(\M) \sim \mu\in\Lambda^2(\V)$
given by $\mu(v_1,v_2)=F(v_1\x v_2)$.
\item $(u,p)\in C^\infty(\V\x\K)\sim (\omega,\mu)\in \Lambda^3(\K\x \V)$
given by $\omega(v_1,v_2,v_3)=p(v_1\x v_2\cdot v_3)$,
$\mu(v_1,v_2,v_3)=u(v_1\x v_2\cdot v_3)$.
\end{itemize}
Under these identifications the we find that
\begin{itemize}
\item $d_0:\Lambda^0(\K)\to \Lambda^1(\K)$ corresponds to the row-wise
gradient $C^\infty(\V)\to C^\infty(\M)$.
\item $S_0:\Lambda^0(\V)\to\Lambda^1(\K)$ corresponds to the inclusion
of $C^\infty(\K) \to C^\infty(\M)$.
\item $d_1\circ S_1^{-1}\circ d_1:\Lambda^1(\K)\to\Lambda^2(\V)$
corresponds to $J=\curl\Xi^{-1}\curl:C^\infty(\M)\to C^\infty(\M)$.
\item $d_2:\Lambda^2(\V)\to\Lambda^3(\V)$ corresponds to the
row-wise divergence $C^\infty(\M)\to C^\infty(\V)$.
\item $S_2:\Lambda^2(\V)\to\Lambda^3(\K)$ corresponds to the
operator $- 2 \skw:C^\infty(\M)\to C^\infty(\K)$.
\end{itemize}
Thus, modulo these identifications and the (unimportant) constant factor
in the last identification, \eqref{elas-seq3} and
\eqref{elas-seq4} are identical. Hence we have established the following 
result.
\begin{thm}
\label{exact-exact}
 When the de~Rham complex \eqref{deRham1} is exact, (i.e., the domain is
contractible), then so is the elasticity complex \eqref{elas-seq3}.
\end{thm}

To end this section, we return to
the operator $S:\Lambda^k(\V)\to\Lambda^{k+1}(\K)$ defined by $S=dK-Kd$.
Let $K':\Lambda^k(\K)\to\Lambda^k(\V)$ be the adjoint of $K$ (with
respect to the Euclidean inner product on $\V$ and the Frobenius
inner product on $\K$), which is
given by $(K'\omega)_x(v_1,\ldots,v_k)=-2\omega_x(v_1,\ldots,v_k)x$.  Define $S':\Lambda^k(\K)\to\Lambda^{k+1}(\V)$
by $S'=dK'-K'd$.
Recall that the wedge product
$\wedge:\Lambda^k\x\Lambda^l\to\Lambda^{k+l}$ is given by
\begin{equation*}
(\omega\wedge\mu)(v_1,\ldots,v_{k+l})=\sum
(\operatorname{sign}\sigma)\omega(v_{\sigma_1},\ldots,v_{\sigma_k})
\mu(v_{\sigma_{k+1}},\ldots,v_{\sigma_{k+l}}),\quad 
\omega\in\Lambda^k, \mu\in\Lambda^l, v_i\in V,
\end{equation*}
where the sum is over the set of all permutations of 
$\{1,\ldots,k+l\}$, for which $\sigma_1<\sigma_2<\cdots<\sigma_k$
and $\sigma_{k+1}<\sigma_{j+2}<\cdots<\sigma_{k+l}$.
This extends as well to differential forms with values in an
inner product space, using the inner product to multiply the terms
inside the summation.  Using the Leibniz rule \eqref{leibniz}, we have
\begin{equation}\label{adj}
(S\omega)\wedge\mu= (-1)^{k}\omega\wedge S'\mu, \quad
\omega\in\Lambda^k(\V),~\mu\in\Lambda^l(\K).
\end{equation}
We thus have
\begin{equation*}
dK\omega\wedge \mu = (-1)^{k+1}K\omega\wedge d\mu +d(K\omega\wedge\mu)
= (-1)^{k+1}\omega\wedge  K'd\mu +d(\omega\wedge K'\mu),
\end{equation*}
and
\[
Kd\omega\wedge \mu =d\omega\wedge K'\mu
= (-1)^{k+1}\omega\wedge  dK'\mu +d(\omega\wedge K'\mu),
\]
Subtracting these two expressions gives \eqref{adj}.

For later reference, we note that, analogously to \eqref{defS}, we have
\begin{equation}\label{defS'}
(S'\omega)(v_1,\ldots,v_{k+1}) = -2\sum_{j=1}^{k+1}(-1)^{j+1}
\omega(v_1,\ldots,\hat v_j,\ldots v_{k+1}) v_j,\quad
\omega\in \Lambda^k(\Omega;\K).
\end{equation}

\section{The discrete construction}
\label{sec:BGG_h}
In this section we derive a discrete version of the elasticity
sequence by adapting the construction of the previous section.
To carry out the construction, we will use two
discretizations of the de~Rham sequence.
For $k=0,1,2,3$, let $\Lambda^k_h$ denote a finite-dimensional
space of $H\Lambda^k$ for which
$d\Lambda^k_h\subset \Lambda^{k+1}_h$, and for which
there exist projections $\Pi_h=\Pi^k_h:\Lambda^k\to\Lambda^k_h$ which
make the following diagram commute:
\begin{equation}\label{deRham-form5}
\begin{CD}
\R  @.\,\hookrightarrow\,@. \Lambda^0 @>d>> \Lambda^1 @>d>>
\Lambda^2 @>d>> \Lambda^3 @>>> 0\\
@. @. @VV\Pi_hV @VV\Pi_hV  @VV\Pi_hV  @VV\Pi_hV 
\\
\R  @.\,\hookrightarrow\,@. \Lambda_h^0 @>d>> \Lambda_h^1 @>d>>
\Lambda_h^2 @>d>> \Lambda_h^3 @>>> 0
\end{CD}
\end{equation}
This is simply the diagram \eqref{diagram4} written in the language of
differential forms.  We do not make a specific choice of the
discretization yet, but, as recalled in \S~\ref{sec:prelim}, there
exist many such discrete de Rham complexes based on piecewise
polynomials.  In fact, as explained in \cite{Arnold-Falk-Winther05_1},
for each polynomial degree $r\ge 0$ we may choose $\Lambda^3_h$ to be
the space of all piecewise polynomial $3$-forms with respect to some
simplicial decomposition of $\Omega$, and construct four such diagrams.
We make the assumption that $\P_1(\Omega)\subset\Lambda^0_h$, which is
true in all the cases mentioned.

Let $\tilde\Lambda^k_h$ be a second set of finite dimensional
spaces with corresponding projection operators $\tilde\Pi_h$ enjoying
the same properties, giving us a second discretization of the de~Rham
sequence.  Supposing a compatibility condition between
these two discretizations, which we describe below, we shall construct a discrete
elasticity complex.

We start with the complex
\begin{equation}\label{deRham-form6}
\begin{CD}
\K\x\V  @.\,\hookrightarrow\,@. \Lambda^0_h(\K)\x\tilde\Lambda^0_h(\V) @>d>> 
\cdots @>d>> \Lambda^3_h(\K)\x\tilde\Lambda^3_h(\V) @>>> 0
\end{CD}
\end{equation}
where $\Lambda^k_h(\K)$ denotes the $\K$-valued analogue of $\Lambda^k_h$
and similarly for $\tilde\Lambda^k_h(\V)$.  For brevity, we henceforth
write $\Lambda^k_h(\W)$ for $\Lambda^k_h(\K)\x\tilde\Lambda^k_h(\V)$.
As a discrete analogue of the operator $K$, we define 
$K_h : \tilde\Lambda_h^k(\V) \to \Lambda_h^k(\K)$ by $K_h = \Pi_h K$
where $\Pi_h$ is the interpolation operator onto $\Lambda_h^k(\K)$.

Next define $S_h =S_{k,h}: \tilde\Lambda_h^k(\V) \to \Lambda_h^{k+1}(\K)$
by $S_h = d K_h - K_hd$, for $k=0,1,2$.  Observe that the discrete version 
of \eqref{dS=-Sd},
\begin{equation}\label{dS=-Sd_h}
dS_h =- S_hd,
\end{equation}
follows exactly as in the continuous case.
From the commutative diagram \eqref{deRham-form5}, we see that
\[
S_h = d\Pi_h K - \Pi_hKd = \Pi_h(dK - Kd) = \Pi_h S.
\] 
Continuing to mimic the continuous case, we define
the automorphism $\Phi_h$ on $\Lambda^k_h(\W)$ by
\[
\Phi_h(\omega,\mu) = (\omega + K_h\mu, \mu),
\]
and the operator $\A_h : \Lambda_h^k(\W) \to \Lambda_h^{k+1}(\W)$
by $\A_h = \Phi_h d \Phi_h^{-1}$, which leads to
\[ 
\A_h(\omega,\mu) = (d\omega - S_h \mu, d\mu).
\]
Thanks to the assumption that $\P_1\subset\Lambda^0_h$, we have
$\Phi_h(\W)=\Phi(\W)$.  Hence,
inserting the isomorphisms $\Phi_h$ into \eqref{deRham-form6}, we obtain
\begin{equation}\label{deRham-form7}
\begin{CD}
\Phi(\W) @.\,\hookrightarrow\,@. \Lambda_h^0(\W) @>\A_h>> 
\Lambda_h^1 (\W)@>\A_h>> \Lambda_h^2(\W) @>\A_h>> \Lambda_h^3(\W)@>>> 0.
\end{CD}
\end{equation}
In analogy to the continuous case, we define
%\begin{equation}\label{discreteE}
%\begin{CD}
%\Phi(\W)  @.\,\hookrightarrow\,@. \Lambda_h^0(\W) @>\A_h>> \Gamma_h^1@>\A_h>>
%\Gamma_h^2@>\A_h>> \Lambda_h^3(\W)@>>> 0.
%\end{CD}
%\end{equation}
\[
\Gamma^1_h = \{\,(\omega, \mu) \in \Lambda^1_h(\W)\,|\, d\omega = S_{1,h}\mu \,\}
,
\quad
\Gamma^2_h = \{\,(\omega, \mu) \in \Lambda^2_h(\W)\,|\, \omega = 0\,\}.
\]
As in the continuous case, we can identify $\Gamma_h^2$ with
$\tilde\Lambda^2_h(\V)$, but, unlike in the continuous case, we cannot
identify $\Gamma^1_h$ with $\Lambda^1_h(\K)$, since we do not require that
$S_{1,h}$ be invertible (it is not in the applications).  Hence, in order to
derive the analogue of the diagram \eqref{projection1} we require a {\bf
surjectivity assumption:}
\begin{equation}\label{boxcond}
\text{\it The operator $S_{1,h}$ maps
$\tilde\Lambda_h^1(\V)$ onto $\Lambda_h^2(\K)$.}
\end{equation}
Under this assumption, the operator $S_h=S_{1,h}$ has a right inverse
$S_h^\dagger$ mapping $\Lambda_h^2(\K)$ into $\Lambda_h^1(\V)$.  This allows
us to define discrete counterparts of the projection operators $\pi^1$ and
$\pi^2$ by
\[
\pi_h^1(\omega, \mu) = (\omega, \mu - S_h^{\dagger}S_h \mu + 
S_h^{\dagger}d \omega), \quad
\pi_h^2(\omega, \mu) = (0 , \mu + d S_h^{\dagger} \omega),
\]
and obtain the discrete analogue of \eqref{projection1}:
\begin{equation}\label{projection2}
\begin{CD}
\Phi(\W)  @.\,\hookrightarrow\,@. \Lambda_h^0(\W) @>\A_h>> \Lambda_h^1 (\W)@>\A_h>>
\Lambda_h^2(\W) @>\A_h>> \Lambda_h^3(\W)@>>> 0\\
@. @. @VVidV @VV\pi_h^1V @VV\pi_h^2V  @VVidV\\
\Phi(\W)  @.\,\hookrightarrow\,@. \Lambda_h^0(\W) @>\A_h>> \Gamma_h^1@>\A_h>>
\Gamma_h^2@>\A_h>> \Lambda_h^3(\W)@>>> 0
\end{CD}
\end{equation}
It is straightforward to check that this diagram commutes.
For example, if $(\omega,\mu) \in \Lambda_h^0(\W)$, then
\begin{align*}
\pi_h^1 \A_h(\omega,\mu) &= \pi_h^1 (d\omega - S_h \mu ,d\mu)
= (d\omega - S_h \mu, d\mu - S_h^{\dagger}S_h d\mu +
S_h^{\dagger}d [d\omega-S_h \mu])
\\
&=  (d\omega - S_h \mu, d\mu - S_h^{\dagger}[S_h d\mu + d S_h \mu])
= \A_h(\omega,\mu),
\end{align*}
where the last equality follows from \eqref{dS=-Sd_h}.  Thus the bottom
row of \eqref{projection2} is a subcomplex of the top row, and
the vertical maps are commuting projections.
In particular, when the top row is exact, so is the bottom.
Thus we have established the following result.
\begin{thm}
  For $k=0, \ldots, 3$, let $\Lambda^k_h$ be a finite dimensional subspace of
  $H \Lambda^k$ for which $d \Lambda^k_h \subset \Lambda^{k+1}_h$ and for
  which there exist projections $\Pi_h = \Pi^k_h: \Lambda^k \to \Lambda^k_h$
  that make the diagram \eqref{deRham-form5} commute.  Let $\tilde
  \Lambda^k_h$ be a second set of finite dimensional spaces with corresponding
  projection operators $\tilde \Pi^k_h$ with the same properties.  If $S_{1,h}
  := d \Pi^1_h K - \Pi^2_h K d$ maps $\tilde \Lambda^1_h(\V)$ onto
  $\Lambda^2_h(\K)$, and the bottom row of \eqref{deRham-form5} is exact for
  both sequences $\Lambda^k_h$ and $\tilde \Lambda^k_h$, then the discrete
  elasticity sequence given by the bottom row of \eqref{projection2} is also
  exact.
\end{thm}

The exactness of the bottom row of \eqref{projection2} suggests that
the following choice of finite element spaces will lead to a stable
discretization of \eqref{emixed-system-h}:
\begin{equation*}
\Sigma_h \sim \tilde\Lambda^2_h(\V), \quad V_h \sim
\tilde\Lambda^3_h(\V), \quad Q_h \sim \Lambda^3_h(\K).
\end{equation*}
In the next section we will make specific choices for the discrete
de~Rham complexes, and then verify the stability in the following section.

For use in the next section, we state the following result, giving
a sufficient condition for the key requirement that $S_{1,h}$ be surjective.
\begin{prop}  If the diagram 
\begin{equation}\label{keydiag}
\begin{CD}
\Lambda^1(\V) @>S_1>> \Lambda^2(\K) \\
@V\tilde\Pi_h^1VV  @V\Pi_h^2VV \\
\tilde\Lambda_h^1(\V) @>S_{1,h}>> \Lambda_h^2(\K)
\end{CD}
\end{equation}
commutes, then $S_{1,h}$ is surjective.
\end{prop}

\section{A family of discrete elasticity complexes}\label{sec:examples}

In this section, we present a family of examples of the general discrete
construction presented in the previous section, by choosing specific
discrete de~Rham complexes. These furnish a family of discrete
elasticity complexes, indexed by an integer degree $r\ge0$. In the next
section we use these complexes to derive finite elements for elasticity.
In the lowest order case, the method will require only piecewise linear
functions to approximate stresses and piecewise constants to
approximate the displacements and multipliers.

We begin by recalling the two principal families of piecewise polynomial
spaces of differential forms, following the presentation in
\cite{Arnold-Falk-Winther05_1}.  We henceforth assume that the domain $\Omega$
is a contractible polyhedron.  Let $\Th$ be a triangulation of $\Omega$. Let
$\Th$ be a triangulation of $\Omega$ by tetrahedra, and set
\begin{align*}
\P_r\Lambda^k(\Th)&=\{\,\omega\in H\Lambda^k(\Omega)\,\,|\,\, \omega_{|_T}\in
\P_r\Lambda^k(T)\quad\forall T\in\Th\,\}, \\
\P_r^+\Lambda^k(\Th)&=\{\,\omega\in H\Lambda^k(\Omega)\,\,|\,\,
\omega_{|_T}\in \P_r^+\Lambda^k(T)\quad\forall T\in\Th\,\}.
\end{align*}
Here $\P_r^+\Lambda^k(T):=\P_r\Lambda^k(T)+\kappa\P_r\Lambda^{k+1}(T)$
where $\kappa:\Lambda^{k+1}(T)\to\Lambda^k(T)$ is the 
\emph{Koszul differential} defined by
\begin{equation*}
(\kappa \omega)_x(v^1, \cdots, v^k) = \omega_x(x,v^1, \cdots, v^k).
\end{equation*}
The spaces $\P_r^+\Lambda^0(\Th)=\P_{r+1}\Lambda^0(\Th)$ correspond to the
usual degree $r+1$ Lagrange piecewise polynomial subspaces of $H^1$, and the
spaces $\P_r^+\Lambda^3(\Th)=\P_r\Lambda^3(\Th)$ correspond to the usual
degree $r$ subspace of discontinuous piecewise polynomials in $L^2(\Omega)$.
For $k=1$ and $2$, the spaces $\P_r^+\Lambda^k(\Th)$ correspond to the
discretizations of $H(\operatorname{curl})$ and $H(\operatorname{div})$,
respectively, presented by N\'ed\'elec in \cite{Nedelec1}, and the spaces
$\P_r\Lambda^k(\Th)$ are the ones presented by N\'ed\'elec in \cite{Nedelec2}.
An element $\omega\in\P_r\Lambda^k(\Th)$ is uniquely determined by the
following quantities:
\begin{equation}\label{dof1}
\int_f \omega\wedge\zeta, \quad \zeta\in \P_{r-d-1+k}^+\Lambda^{d-k}(f),
\quad f\in \Delta_d(\Th), \quad k\le d \le 3.
\end{equation}
Here $\Delta_d(\Th)$ is the set of vertices, edges, faces, or tetrahedra
in the mesh $\Th$, for $d=0,1,2,3$, respectively, and
for $r<0$, we interpret $\P_r^+\Lambda^k(T)=\P_r\Lambda^k(T)=0$.
Note that for $\omega\in\Lambda^k$, $\omega$ naturally restricts
on the face $f$ to an element of $\Lambda^k(f)$.  Therefore, for
$\zeta\in\Lambda^{d-k}(f)$, the wedge product
$\omega\wedge\zeta$ belongs to $\Lambda^d(f)$ and hence the integral of
$\omega\wedge\zeta$ on the $d$-dimensional face $f$ of $T$ is a well-defined
and natural quantity.  Using the quantities in \eqref{dof1}, we
obtain a projection operator from $\Lambda^k$ to $\P_r\Lambda^k(\Th)$.

Similarly, an element $\omega\in\P_r^+\Lambda^k(\Th)$ is uniquely determined by
\begin{equation}\label{dof2}
\int_f \omega\wedge\zeta, \quad \zeta\in \P_{r-d+k}\Lambda^{d-k}(f),
\quad f\in \Delta_d(\Th), \quad k\le d \le 3,
\end{equation}
and so these quantities determine a projection.

If $\X$ is a vector space, we use the notation $\P_r\Lambda^k(\Th;\X)$ and
$\P^+_r\Lambda^k(\Th;\X)$ to denote the corresponding spaces of piecewise
polynomial differential forms with values in $\X$. Furthermore, if $\X$ is an
inner product space, the corresponding degrees of freedom are given by
\eqref{dof1} and \eqref{dof2}, but where the test spaces are replaced by the
corresponding $\X$ valued spaces.

To carry out the construction described in the previous section
we need to choose the two sets of spaces $\Lambda^k_h$ and $\tilde\Lambda^k_h$
for $k=0,1,2,3$.  We fix $r\ge0$ and set $\Lambda^k_h=\P_r^+\Lambda^k(\Th)$,
$k=0,1,2,3$, and $\tilde\Lambda^0_h=\P_{r+2}\Lambda^0(\Th)$,
$\tilde\Lambda^1_h=\P_{r+1}^+\Lambda^1(\Th)$, $\tilde\Lambda^2_h=
\P_{r+1}\Lambda^2(\Th)$, and $\tilde\Lambda^3_h=P_r\Lambda^3_h$.
As explained in \cite{Arnold-Falk-Winther05_1},
both these choices give a discrete de~Rham sequence with commuting projections,
i.e., a diagram like \eqref{deRham-form5} makes sense and is commutative.

We establish the key surjectivity assumption for our choice
of spaces, by verifying the commutativity of \eqref{keydiag}.
\begin{lem}
\label{Rhidentproof}
Let $\tilde\Lambda^1_h(\V)=\P_{r+1}^+\Lambda^1(\Th;\V)$ and
$\Lambda^2_h(\K)
=\P_r^+\Lambda^2(\Th;\K)$
with projections $\tilde\Pi^1_h$,
$\Pi^2_h$ defined via the corresponding vector valued moments of the
form
\eqref{dof1} and \eqref{dof2}. If
$S_{1,h}=\Pi^2_h S_1$ then 
\begin{equation}\label{Rhident}
S_{1,h}\tilde\Pi_h^1=\Pi_h^2 S_1,
\end{equation}
and so $S_{1,h}$ is surjective.
\end{lem}
\begin{proof}
We must show that $(\Pi^2_h S_1-S_{1,h}\tilde\Pi^1_h)\sigma=0$ for
$\sigma\in\Lambda^1(\V)$. Defining $\omega=(I-\tilde\Pi^1_h)\sigma$,
the required condition becomes $\Pi^2_hS_1\omega=0$.  Since
$\tilde\Pi^1_h\omega=0$, we have
\begin{equation}\label{dof1v}
\int_f \omega\wedge\zeta = 0 , \quad
\zeta\in\P_{r+2-d}\Lambda^{d-1}(f;\V), \quad
f\in\Delta_d(\Th), \quad 2\le d\le 3,
\end{equation}
(in fact \eqref{dof1v} holds for $d=1$ as well, but
this is not used here).
We must show that \eqref{dof1v} implies
\begin{equation}\label{dof2v}
\int_f S_1\omega\wedge\mu = 0 , \quad
\mu\in\P_{r+2-d}\Lambda^{d-2}(f;\K),\quad
f\in\Delta_d(\Th), \quad 2\le d\le 3.
\end{equation}
From \eqref{adj}, we have
$S_1\omega\wedge\mu= - \omega\wedge\zeta$, where $\zeta=S'_{d-2}\mu
\in \P_{r+2-d}\Lambda^{d-1}(f;\V)$ for
$\mu\in\P_{r+2-d}\Lambda^{d-2}(f;\K)$, as is evident from \eqref{defS'}.
Hence \eqref{dof2v} follows from \eqref{dof1v}.
\end{proof}

\section{Stable mixed finite elements for elasticity}

Based on the discrete elasticity complexes just constructed, we obtain
mixed finite element spaces for the formulation \eqref{emixed-system-h}
of the elasticity equations by choosing $\Sigma_h$, $V_h$, and $Q_h$
as the spaces of matrix and vector fields corresponding to appropriate
spaces of forms in the $\K$- and $\V$-valued de~Rham sequences used
in the construction.  Specifically, these are
\begin{equation}\label{spaces}
\Sigma_h \sim \P_{r+1}\Lambda^2(\Th;\V), \quad
V_h \sim \P_r\Lambda^3(\Th;\V), \quad
Q_h \sim \P_r\Lambda^3(\Th;\K).
\end{equation}

In other terminology, $\Sigma_h$ may be thought of
as the product of three copies of the
N\'ed\'elec $H(\operatorname{div})$ space of the second kind of degree
$r+1$, and $V_h$ and $Q_h$ are spaces of all piecewise polynomials of
degree at most $r$ with values in $\K$ and $\V$, respectively, with
no imposed interelement continuity.  In this section, we establish 
stability and convergence for this finite element method.

The stability of the method requires the two conditions (A$1$) and (A$2$)
stated in \S~2.  The first condition is obvious since, by construction,
$\div\Sigma_h\subset V_h$, i.e., $d\P_{r+1}\Lambda^2(\Th;\V)\subset
\P_r\Lambda^3(\Th;\V)$.
The condition (A$2$) is more subtle.  We will prove a stronger version, namely
\begin{itemize}
\setlength{\itemindent}{20pt}
\item[\quad(A$2'$)] for all nonzero $(v,q) \in V_h \x Q_h$,
there exists nonzero $\tau\in\Sigma_h$ with $\div\tau=v$,
$2\Pi_{Q_h}\skw\tau = q$ and
\begin{equation*}
\|\tau\|_{\div} \le c (\|v\|+\|q\|),
\end{equation*}
where $\Pi_{Q_h}$ is the $L^2$ projection into $Q_h$ and $c$ is a
constant.
\end{itemize}
Recalling that $\Gamma^2_h=0\times \P_{r+1}\Lambda^2(\Th;\V)$
and $\A_h(0, \sigma) = (-S_{2,h} \sigma, d \sigma)$, and that
the operator $S_2$ corresponds on the matrix level
to $-2\skw$, we restate
condition (A$2'$) in the language of differential forms.
\begin{thm}
\label{thm1-A2}
Given $(\omega, \mu) \in \P_r\Lambda^3(\Th;\K) \times \P_r\Lambda^3(\Th;\V)$,
there exists $\sigma \in \P_{r+1}\Lambda^2(\Th;\V)$ such that
$\A_h(0, \sigma) = (\omega, \mu)$ and
\begin{equation}\label{bnd}
\|\sigma\|_{H\Lambda} \le c (\|\omega\| + \|\mu\|),
\end{equation}
where the constant $c$ is independent of $\omega,\mu$ and $h$.
\end{thm}

Before proceeding to the proof, we need to establish some bounds
on projection operators.  We do this for the corresponding scalar-valued
spaces.  The extensions to vector-valued spaces are straightforward.
First we claim that
\begin{equation}\label{projbnds}
\|\tilde \Pi^2_h\eta\|\le c\|\eta\|_1 \quad\forall\eta\in H^1\Lambda^2,\quad
\|\Pi^3_h\omega\|\le c\|\omega\|_0\quad\forall\omega\in H^1\Lambda^3.
\end{equation}
Here the constant may depend on the shape regularity of the mesh, but
not on the meshsize. The second bound is obvious (with $c=1$), since
$\Pi^3_h$ is just the $L^2$ projection.  The first bound follows by
a standard scaling argument.  Namely, let $\hat T$ denote the reference
simplex.  For any $\hat\beta\in \P_{r+1}\Lambda^2(\hat T)$, we have
\begin{equation}\label{b1}
\|\hat\beta\|_{0,\hat T} \le c (\sum_{\hat f}\sum_{\hat\mu}|\int_{\hat
f} \hat\beta\wedge\hat\mu| +\sum_{\hat\zeta}|\int_{\hat
T}\hat\beta\wedge\hat\zeta|),
\end{equation}
where $\hat f$ ranges over the faces of $\hat T$, $\hat\mu$ over a basis
for $\P^+_r\Lambda^0(\hat f)$, and $\hat\zeta$ over a basis for
$\P^+_{r-1}\Lambda^1(\hat T)$.  This is true because the integrals on the
right hand side of \eqref{b1} form a set of degrees of freedom for
$\hat\beta\in \P_{r+1}\Lambda^2(\hat T)$ (see \eqref{dof1}), and so we
may use the equivalence of all norms on this finite dimensional space.
We apply this result with $\hat\beta=\hat\Pi^2_h\hat\eta$,
where
$\hat\Pi^2_h$ is the projection defined to preserve
the integrals on the right hand
side of \eqref{b1}.  It follows that
\begin{equation*}
\|\hat\Pi^2_h\hat\eta\|_{0,\hat T}
\le  c(\sum_{\hat f}\sum_{\hat\mu}|\int_{\hat
f} \hat\eta\wedge\hat\mu| +\sum_{\hat\zeta}|\int_{\hat
T}\hat\eta\wedge\hat\zeta|)\le c\|\hat\eta\|_{1,\hat T},
\end{equation*}
where we have used a standard trace inequality in the last step.
Next, if $T$ is an arbitrary simplex and $\eta\in H^1\Lambda^2(T)$,
we map the reference simplex $\hat T$ onto $T$ by an affine map
$\hat x\mapsto B\hat x +b$, and define $\hat\eta\in H^1\Lambda^2(\hat T)$
by
\begin{equation*}
\hat\eta_{\hat x}(\hat v_1, \hat v_2) = \eta_x(B\hat v_1, B \hat v_2),
\end{equation*}
for any $x=B\hat x+b\in T$ and any vectors $\hat v_1, \hat v_2$.  It is easy
to check that $\widehat{\tilde\Pi^2_h\eta}=\hat\Pi^2_h\hat\eta$, and that
\begin{equation*}
\|\tilde\Pi^2_h\eta\|_{0,T}\le c \|\hat\Pi^2_h\hat\eta\|_{0,\hat T}
\le c\|\hat\eta\|_{1,\hat T}\le c (\|\eta\|_{0,T} + h|\eta|_{1,T})
\le c\|\eta\|_{1,T}.
\end{equation*}
Squaring and adding this over all the simplices in the mesh $\Th$ gives the
first bound in \eqref{projbnds}.

We also need a bound on the projection of a form in $H^1\Lambda^1$ into
$\tilde\Lambda^1_h=\P^+_{r+1}\Lambda^1(\Th)$.  However, the projection
operator $\tilde\Pi^1_h$ is \emph{not} bounded on $H^1$,  because its
definition involves integrals over edges. A similar problem has arisen
before (see, e.g., \cite{Arnold-Winther02}), and we use the same
remedy.  Namely we start by defining an operator
$\tilde\Pi^1_{0h}:H^1\Lambda^1\to \P^+_{r+1}\Lambda^1(\Th)$ by
the conditions
\begin{align}\label{Pi01}
\int_T\tilde\Pi^1_{0h}\omega\wedge\zeta&=\int_T\omega\wedge\zeta, \quad
\zeta\in\P_{r-1}\Lambda^2(T),\quad T\in\Th,
\\\label{Pi02}
\int_f\tilde\Pi^1_{0h}\omega\wedge\zeta&=\int_f\omega\wedge\zeta, \quad
\zeta\in\P_r\Lambda^1(f),\quad f\in\Delta_2(\Th),
\\\label{Pi03}
\int_e\tilde\Pi^1_{0h}\omega\wedge\zeta&=0, \quad
\zeta\in\P_{r+1}\Lambda^0(e), \quad e\in\Delta_1(\Th).
\end{align}
Note that, in contrast to $\tilde\Pi^1_h$, in the definition
of $\tilde\Pi^1_{0h}$, we have set the troublesome edge degrees of freedom to
zero.
Let $\hat\Pi^1_0:H^1\Lambda^1(\hat T)\to \P^+_{r+1}\Lambda^1(\hat T)$
be defined analogously on the reference element.

Now for $\hat\rho\in
H^1\Lambda^1(\hat T)$, $\hat d\hat\Pi^1_0\hat\rho\in\P_{r+1}\Lambda^2(\hat
T)$, so
\begin{equation*}
\|\hat d\hat\Pi^1_0\hat\rho\|_{0,\hat T}
\le c (\sum_{\hat f}\sum_{\hat\mu}
  |\int_{\hat f} \hat d\hat\Pi^1_0\hat\rho\wedge\hat\mu|
 +\sum_{\hat\zeta}|\int_{\hat T}\hat d\hat\Pi^1_0\hat\rho\wedge\hat\zeta|),
\end{equation*}
where again $\hat f$ ranges over the faces of $\hat T$, $\mu$ over a
basis of $\P^+_r\Lambda^0(\hat f)$, and
$\zeta$ over a basis of $\P^+_{r-1}\Lambda^1(\hat T)$.
But
\begin{equation*}
\int_{\hat f} \hat d\hat\Pi^1_0\hat\rho\wedge\hat\mu
=\int_{\hat f}\hat\Pi^1_0\hat\rho \wedge \hat d\hat\mu
=\int_{\hat f}\hat\rho \wedge \hat d\hat\mu,
\end{equation*}
where we have used Stokes theorem and the fact that the vanishing of
the edge quantities in the definition of $\hat\Pi^1_0$ to obtain the
first equality, and the face degrees of freedom entering the definition
of $\hat\Pi^1_0$ to obtain the second.  Similarly,
\begin{equation*}
\int_{\hat T}\hat d\hat\Pi^1_0\hat\rho\wedge\hat\zeta
=\int_{\hat T}\hat\Pi^1_0\hat\rho\wedge\hat d\hat\zeta
+\int_{\partial\hat T}\hat\Pi^1_0\hat\rho\wedge\hat\zeta
=\int_{\hat T}\hat\rho\wedge\hat d\hat\zeta
 +\int_{\partial\hat T}\hat\rho\wedge\hat\zeta
 =\int_{\hat T}\hat d\hat\rho\wedge\hat\zeta.
\end{equation*}
It follows that
\begin{equation*}
\|\hat d\Pi^1_0\hat\rho\|_{0,\hat T} \le c \|\hat\rho\|_{1,\hat T},
\quad \rho\in H^1\Lambda^1(\hat T).
\end{equation*}
When we scale this result to an arbitrary simplex and add over the mesh,
we obtain
\begin{equation*}
\|d\Pi^1_{0h}\rho\| \le c (h^{-1}\|\rho\| + \|\rho\|_1),\quad
\rho\in H^1\Lambda^1(\Omega).
\end{equation*}

To remove the problematic $h^{-1}$ in the last estimate, we introduce
the Clement interpolant $R_h$ mapping $H^1\Lambda^1$ into continuous
piecewise linear $1$-forms (still following \cite{Arnold-Winther02}).  Then 
\begin{equation*}
\|\rho-R_h\rho\|+h\|\rho-R_h\rho\|_1 \le c h \|\rho\|_1, \quad
\rho \in H^1\Lambda^1.
\end{equation*}
Defining $\bar\Pi^1_h:H^1\Lambda^1\to \P^+_{r+1}\Lambda^1_h$ by
\begin{equation}\label{Pibar}
\bar\Pi^1_h = \tilde\Pi^1_{0h}(I-R_h)+ R_h,
\end{equation}
we obtain
\begin{equation*}
\|d\bar\Pi^1_h\rho\|\le \|d\tilde\Pi^1_{0h}(I-R_h)\rho\| + \|dR_h\rho\|
\le  c (h^{-1}\|(I-R_h)\rho\| + \|(I-R_h)\rho\|_1+ \|dR_h\rho\|)
\le c\|\rho\|_1.
\end{equation*}
Thus we have shown that
\begin{equation}\label{p1bnd}
\|d\bar\Pi^1_h\rho\|\le c\|\rho\|_1, \quad \rho\in H^1\Lambda^1.
\end{equation}

Having modified $\tilde\Pi^1_h$ to obtain the bounded operator
$\bar\Pi^1_h$, we now verify that the key property \eqref{Rhident}
in Lemma~\ref{Rhidentproof}
carries over to
\begin{equation}\label{Rhidentb}
S_{1,h}\bar\Pi_h^1=\Pi_h^2 S_1,
\end{equation}
where we now use the vector-valued forms of the projection operators.
It follows
easily from \eqref{Pibar}, \eqref{Pi01}, and \eqref{Pi02} that
\eqref{dof1v} holds with $\omega=(I-\bar\Pi^1_h)\sigma$, so
that the proof of \eqref{Rhidentb} is the same as for \eqref{Rhident}.

We can now give the proof of Theorem \ref{thm1-A2}.
\begin{proof}[Proof of Theorem \ref{thm1-A2}]
Given $\mu\in \P_r\Lambda^3(\Th;\V)$ there exists $\eta\in H^1\Lambda^2(\V)$
such that $d\eta=\mu$ with the bound $\|\eta\|_1\le c\|\mu\|$
(since $d$ maps $H^1\Lambda^2$ onto
$L^2\Lambda^3$).  Similarly, given $\omega\in \P_r\Lambda^3(\Th;\K)$ there
exists $\tau\in H^1\Lambda^2(\K)$ with $d\tau=\omega+S_{2,h}\tilde\Pi^2_h\eta$
with the bound $\|\tau\|_1\le c\|\omega+S_{2,h}\tilde\Pi^2_h\eta\|$.
Let $\rho=S_1^{-1}\tau$ (recall that $S_1$ is an isomorphism) and set
\[
\sigma = d\bar\Pi^1_h\rho +\tilde\Pi^2_h\eta .
\]
We will now show that $\A_h(0, \sigma) = (\omega, \mu)$.
From the definition of $\sigma$, we have
\[
-S_{2,h}\sigma =-S_{2,h} d\bar\Pi^1_h\rho -S_{2,h}\tilde\Pi^2_h\eta .
\]
Then, using \eqref{dS=-Sd_h}, \eqref{Rhidentb}, and the commutativity
$d\Pi_h=\Pi_hd$, we see
\begin{multline*}
S_{2,h} d\bar\Pi^1_h\rho =-d S_{1,h}\bar\Pi^1_h\rho
=-d \Pi^2_h S_1\rho
\\
=-d\Pi^2_h\tau=-\Pi^3_h d\tau =
-\Pi^3_h(\omega+S_{2,h}\tilde\Pi^2_h\eta)= -\omega-S_{2,h}\tilde\Pi^2_h\eta.
\end{multline*}
Combining, we get $-S_{2,h}\sigma=\omega$ as desired.
Further, from the commutativity $d\tilde\Pi_h=\tilde\Pi_hd$ and the
definition of $\eta$, we get
\[
d\sigma=d\tilde\Pi^2_h\eta=\tilde\Pi^3_hd\eta=\tilde\Pi^3_h\mu=\mu,
\]
and so we have established that $\A_h(0,\sigma)=(\mu,\omega)$.

It remains to prove the bound \eqref{bnd}.  Using
\eqref{projbnds}, we have
\begin{equation*}
\|S_{2,h}\tilde\Pi^2_h\eta\|=\|\Pi^3_h S_2\tilde\Pi^2_h\eta\|
\le c \|S_2\tilde\Pi^2_h\eta\| \le c\|\tilde\Pi^2_h\eta\|
\le c\|\eta\|_1\le c\|\mu\|.
\end{equation*}
Thus $\|\rho\|_1\le c\|\tau\|_1\le c(\|\omega\|+\|\mu\|)$.
Using \eqref{p1bnd}, we then get
$\|d\bar\Pi^1_h\rho\|\le c\|\rho\|_1\le c(\|\omega\|+\|\mu\|)$,
and, using \eqref{projbnds}, that $\|\tilde\Pi^2_h\eta \|\le c\|\eta\|_1\le
c\|\mu\|$.  Therefore $\|\sigma\|\le c(\|\omega\|+\|\mu\|)$,
while $\|d\sigma\|=\|\mu\|$, and thus we have the desired bound
\eqref{bnd}.
\end{proof}

We have thus verified the stability conditions (A$1$) and (A$2$), and so may
apply the standard theory of mixed methods (cf.~\cite{Brezzi},
\cite{Brezzi-Fortin}, \cite{douglas-roberts}, \cite{falk-osborn}) and
standard results about approximation by finite element spaces to obtain
convergence and error estimates.
\begin{thm}
\label{convergence-thm}
Suppose $(\sigma, u, p)$ is the solution of the elasticity
system \eqref{emixed-system}
and $(\sigma_h, u_h, p_h)$ is the solution of discrete
system \eqref{emixed-system-h},
where the finite element spaces $\Sigma_h$, $V_h$, and $Q_h$
are given by \eqref{spaces} for some integer $r\ge 0$. Then there is a
constant $C$, independent of $h$, such that
\begin{gather*}
\|\sigma - \sigma_h\|_{\div} + \|u-u_h\| + \|p-p_h\| \le
C \inf_{\tau_h \in \Sigma_h, v_h \in V_h, q_h \in Q_h}
(\|\sigma - \tau_h\|_{\div} + \|u-v_h\| + \|p-q_h\|),
\\
\|\sigma- \sigma_h\| + \|p-p_h\| + \|u_h - \tilde \Pi^n_h u\|
 \le C (\|\sigma- \tilde \Pi^{n-1}_h \sigma\| + \|p - \Pi^n_h p\|),
\\
\|u-u_h\| \le C(\|\sigma- \tilde \Pi^{n-1}_h \sigma\| 
+ \|p - \Pi^n_h p\| + \|u- \tilde \Pi^n_h u\|),
\\
\|d (\sigma - \sigma_h)\| = \|d \sigma - \tilde \Pi^n d \sigma\|.
\end{gather*}
If $u$ and $\sigma$ are sufficiently smooth, then
\begin{equation*}
\|\sigma - \sigma_h\| + \|u-u_h\| + \|p-p_h\| \le Ch^{r+1} \|u\|_{r+2},
\qquad \|\div(\sigma - \sigma_h)\| \le Ch^{r+1} \|\div \sigma\|_{r+1}.
\end{equation*}
\end{thm}

\begin{remark} Note that the errors
$\|\sigma- \sigma_h\|$ and $\|u_h - \tilde \Pi^n_h u\|$ depend on the
approximation of both $\sigma$ and $p$.  For the choices made here,
the approximation of $p$ is one order less than the approximation
of $\sigma$, and thus we do not obtain improved estimates, as one
does in the case of the approximation of Poisson's equation, where
the extra variable $p$ does not enter.
\end{remark}

\section{A simplified element}\label{sec:simplified}
Recall that the lowest order element in the stable family described 
above, for a discretization based on 
\eqref{emixed-system}, is of the form
\[
\Sigma_h \sim \P_1 \Lambda^2(\Th;\V), \qquad V_h \sim
\P_0 \Lambda^3(\Th;\V),
\qquad Q_h \sim \P_0 \Lambda^3(\Th;\K).
\]
The purpose of this section 
is to present a stable element which is slightly simpler than this one.
More precisely, the spaces $V_h$ and $Q_h$ are unchanged, but
$\Sigma_h$
will be simplified from full linears to matrix fields whose 
tangential--normal components on each two dimensional face of a
tetrahedron are only a reduced space of linears.

We will still adopt the notation of differential forms. By examining
the proof of Theorem~\ref{thm1-A2},  we realize that we do not use the
complete sequence \eqref{deRham-form6} for the given spaces. We
only use the sequences
\begin{equation}\label{reduced-seq}
\begin{array}{rl}
\begin{CD}
&\P_0^+ \Lambda^2(\Th;\K) @>d>> \P_0 \Lambda^3(\Th;\K)@>>> 0,
\end{CD}\\
\begin{CD}
\P_{1}^+\Lambda^1 (\Th;\V)@>d>>
&\P_{1} \Lambda^2(\Th;\V) @>d>> \P_0 \Lambda^3(\Th;\V)@>>> 0.
\end{CD}
\end{array}
\end{equation}
The purpose here is to show that it 
is possible to choose subspaces of some of the spaces in \eqref{reduced-seq}
such that the desired properties still hold. More precisely,
compared to \eqref{reduced-seq}, the spaces $\P_{1}^+\Lambda^1 (\Th;\V)$
and $\P_{1} \Lambda^2(\Th;\V)$ are simplified, while the three
others remain unchanged.  If we denote by $\P_{1,-}^+\Lambda^1 (\Th;\V)$
and $\P_{1,-} \Lambda^2(\Th;\V)$ the simplifications of the spaces
$\P_{1}^+\Lambda^1 (\Th;\V)$ and $\P_{1} \Lambda^2(\Th;\V)$, respectively,
then the properties we need are that:
\begin{equation}\label{simpexact}
\begin{CD}
\P_{1,-}^+\Lambda^1 (\Th;\V)@>d>>
&\P_{1,-} \Lambda^2(\Th;\V) @>d>> \P_0 \Lambda^3(\Th;\V)@>>> 0
\end{CD}
\end{equation}
is a complex and that the surjectivity assumption \eqref{boxcond} holds, i.e., 
$S_h = S_{1,h}$ maps the space $\P_{1,-}^+\Lambda^1 (\Th;\V)$ onto $\P_0^+
\Lambda^2(\Th;\K)$.  Note that if $\P_{0}^{+} \Lambda^2(\Th;\V) \subset
\P_{1,-} \Lambda^2(\Th;\V)$, then $d$ maps $\P_{1,-} \Lambda^2(\Th;\V)$
onto $\P_0 \Lambda^3(\Th;\V)$.

We first show how $\P_{1,-}^+\Lambda^1 (\Th;\V)$ can be
constructed as a subspace of $\P_{1}^+\Lambda^1 (\Th;\V)$.  Since the
construction is done locally on each tetrahedron, we will show how to
construct a space $\P_{1,-}^+\Lambda^1 (T;\V)$ as a subspace of
$\P_{1}^+\Lambda^1 (T;\V)$.  We begin by recalling that
the face degrees of freedom of $\P_{1}^+\Lambda^1 (T;\V)$ have the form
\begin{equation*}
\int_f \omega \wedge \mu, 
\qquad \mu \in \P_0 \Lambda^1(f, \V).
\end{equation*}
We then observe that this six dimensional space can be decomposed into
\[
\P_0 \Lambda^1(f; \V) = \P_0 \Lambda^1(f; T_f) + \P_0 \Lambda^1(f; N_f),
\]
i.e. into forms with values in the tangent space to $f$, $T_f$ or the 
normal space $N_f$. This is a $4+2$ dimensional decomposition.
Furthermore, 
\[
\P_0 \Lambda^1(f; T_f) = \P_0 \Lambda_{\sym}^1(f; T_f)
+ \P_0 \Lambda_{\skw}^1(f; T_f),
\]
where $\mu \in \P_0 \Lambda^1(f; T_f)$ is in $\P_0 \Lambda_{\sym}^1(f;
T_f)$
if and only if $\mu(s)\cdot t = \mu(t)\cdot s$ for orthonormal tangent
vectors
$s$ and $t$. Finally, we obtain a $3+3$ dimensional decomposition
\[
\P_0 \Lambda^1(f; \V) = \P_0 \Lambda_{\sym}^1(f;T_f)
+ \P_0 \Lambda_{\skw}^1(f; \V),
\]
where 
\[
\P_0 \Lambda_{\skw}^1(f; \V)= \P_0 \Lambda_{\skw}^1(f; T_f)
+ \P_0 \Lambda^1(f; N_f).
\]
In more explicit terms, if $\mu \in \P_0 \Lambda^1(F; \V)$
has the form
\begin{equation*}
\mu(q) = (a_1 t + a_2 s + a_3 n) q \cdot t + (a_4 t + a_5 s + a_6 n) q \cdot s,
\end{equation*}
where $t$ and $s$ are orthonormal tangent vectors on $f$, $n$ is
the unit normal to $f$, and $q$ is a tangent vector on $f$,
then we can write
$\mu = \mu_1 + \mu_2$, with $\mu_1 \in \P_0 \Lambda^1_{\sym}(f;\V)$
and $\mu_2 \in \P_0 \Lambda^1_{\skw}(f;\V)$, where
\begin{align*}
\mu_1(q) &=  \left(a_1t + \frac{a_2+a_4}{2} s\right) q \cdot t 
+  \left(\frac{a_2+a_4}{2} t + a_5 s\right) q \cdot s,
\\
\mu_2(q) &=  \left(\frac{a_2-a_4}{2} s + a_3 n\right) q \cdot t 
+ \left(\frac{a_4-a_2}{2} t + a_6 n\right)
 q \cdot s.
\end{align*}

The reason for this particular decomposition of the degrees of freedom is
that if we examine the proof of
Lemma~\ref{Rhidentproof}, where equation \eqref{Rhident} is established, we
see that the only degrees of freedom that are used for an element $\omega \in
\P_{1}^+\Lambda^1 (T;\V)$ are the subset of the face degrees of freedom
given by:
\begin{equation*}
\int_f \omega \wedge (S_0^{\prime} \nu), 
\qquad \nu \in \P_0 \Lambda^0(f; \K).
\end{equation*} 
However, for $\nu \in \P_0 \Lambda^0(f;\K)$, $\mu = S_0^{\prime} \nu$ is given
by $\mu(q) = \nu q$. Since the general element $\nu \in \P_0 \Lambda^0(\K)$
can be written in the form $b_1(t s^T - s t^T) + b_2 (n t^T - t n^T) + b_3 (n
s^T - s n^T)$, $\nu q = (- b_1 s + b_2 n) q \cdot t
+ (b_1 t + b_3 n) q \cdot s$ for
$q$ a tangent vector, and thus $\mu \in \P_0 \Lambda^1_{\skw}(f; \V)$.
Hence, we have split the degrees of freedom into three on each face that we
need to retain for the proof of Lemma~\ref{Rhidentproof} and three on each
face that are not needed.  The reduced space $\P_{1,-}^+ \Lambda^1 (T;\V)$
that we now construct has two properties.  The first is that it still contains
the space $\P_{1} \Lambda^1 (T;\V)$ and the second is that the unused face
degrees of freedom are eliminated (by setting them equal to zero).  We can
achieve these conditions by first writing an element $\omega \in \P_{1}^+
\Lambda^1 (T;\V)$ as $\omega = \Pi_h \omega + (I- \Pi_h) \omega$, where
$\Pi_h$ denotes the usual projection operator into $\P_{1} \Lambda^1 (T;\V)$
defined by the edge degrees of freedom.  Then the elements in $(I- \Pi_h)
\P_{1}^+ \Lambda^1 (T;\V)$ will satisfy
\begin{equation*}
\int_e \omega \wedge \mu =0, \quad \mu \in \P_1 \Lambda^0(e;\V),
\quad e \in \Delta_1(T),
\end{equation*}
i.e., their traces on the edges will be zero.  Thus, they are completely
defined by the face degrees of freedom
\begin{equation*}
\int_f \omega \wedge \mu, 
\qquad \mu \in \P_0 \Lambda^1(f;\V), \quad f \in \Delta_2(T).
\end{equation*}
Since this is the case, we henceforth denote $(I- \Pi_h) \P_{1}^+ \Lambda^1
(T;\V)$ by $\P_{1,f}^+ \Lambda^1 (T;\V)$.  

We then define our reduced space
\begin{equation*}
\P_{1,-}^+ \Lambda^1 (T;\V)
= \P_{1} \Lambda^1 (T;\V) + \P_{1,f,-}^+ \Lambda^1 (T;\V),
\end{equation*}
where $\P_{1,f,-}^+ \Lambda^1 (T;\V)$ denotes the set of forms
$\omega \in \P_{1,f}^+ \Lambda^1 (T;\V)$ satisfying
\begin{equation*}
\int_f \omega \wedge \mu = 0, 
\qquad \mu \in \P_0 \Lambda^1_{\sym}(f;\V),
\end{equation*}
i.e., we have set the unused degrees of freedom to be zero.

Then
\begin{equation*}
\P_{1,-}^+\Lambda_h^1 (\Th;\V)= \{ \omega \in \P_{1}^+\Lambda^1(\Th;\V)
: \, \omega|_T \in \P_{1,-}^+\Lambda^1(T;\V), \, \forall T \in \Th
\}.
\end{equation*}
The degrees of freedom for this space are then
given by
\begin{equation}
\label{dofs-reduced1n}
\int _e \omega \wedge \mu, \quad \mu \in \P_1 \Lambda^0(e;\V), e \in
\Delta_1(T),
\quad 
\int_f \omega \wedge \mu, \quad \mu \in \P_0 \Lambda_{skw}^1(f;\V),
f \in \Delta_2(T).
\end{equation}
It is clear from this definition that the space $\P_{1,-}^+\Lambda^1
(T;\V)$ will have 48 degrees of freedom (36 edge degrees of freedom and 12
face degrees of freedom).  The unisolvency of this space follows immediately
from the unisolvency of the spaces $\P_{1} \Lambda^1(T;\V)$ and
$\P_{1,f,-}^{+} \Lambda^1(T;\V)$.

The motivation for this choice of the space $\P_{1,-}^+\Lambda_h^1 (\Th;\V)$
is that it easily leads to a definition of the space
$\P_{1,-} \Lambda^2(\Th;\V)$ that satisfies the property that
\eqref{simpexact} is a complex.
We begin by defining
\begin{equation*}
\P_{1,-} \Lambda^2 (T;\V)
= \P_{0}^+ \Lambda^2 (T;\V) + d \P_{1,f,-}^+ \Lambda^1 (T;\V).
\end{equation*}
It is easy to see that this space will have 24 face degrees of freedom.
Note this is a reduction of the space $\P_{1} \Lambda^2 (T;\V)$, since
\begin{equation*}
\P_{1} \Lambda^2 (T;\V)
= \P_{0}^+ \Lambda^2 (T;\V) + d \P_{1,f}^+ \Lambda^1 (T;\V).
\end{equation*}
We then define
\begin{equation*}
\P_{1,-} \Lambda^2 (\Th;\V)= \{ \omega \in \P_{1} \Lambda^2(\Th;\V)
: \, \omega|_T \in \P_{1,-} \Lambda^2(T;\V), \, \forall T \in \Th
\}.
\end{equation*}
It is clear that $\P_{0}^+ \Lambda^2(\Th;\V) \subset \P_{1,-}
\Lambda^2(\Th;\V)$.  The fact that the complex \eqref{simpexact} is exact
now follows directly from the fact that the complex
\begin{equation}
\begin{CD}
\P_{1} \Lambda^1 (T;\V)@>d>>
&\P_{0}^+ \Lambda^2(T;\V) @>d>> \P_0 \Lambda^3(T;\V)@>>> 0
\end{CD}
\end{equation}
is exact and the definition
\begin{equation*}
P_{1,-} \Lambda^1(T;\V) = \P_{0}^+ \Lambda^1(T;\V) 
+ \P_{1,f,-}^+ \Lambda^1(T;\V).
\end{equation*}

We will define appropriate degrees of freedom for the space 
$\P_{1,-} \Lambda^2 (T;\V)$ by using a subset
of the 36 degrees of freedom for $\P_{1} \Lambda^2 (T;\V)$, i.e., of $\int_f
\omega \wedge \mu$, $\mu \in \P_1 \Lambda^0(f;\V)$.  In particular,
we take as degrees of freedom for $\P_{1,-} \Lambda^2 (T;\V)$,
\begin{equation*}
\int_f \omega \wedge \mu,\quad \mu \in \P_{1,skw} \Lambda^0(f;\V),
\qquad \forall f \in \Delta_2(T),
\end{equation*}
where $\P_{1,skw} \Lambda^0(f;\V)$ denotes the set of
$\mu \in \P_1 \Lambda^0(f;\V)$ that satisfy
$d \mu \in \P_0 \Lambda^1_{skw}(f;\V)$. It is easy to check that
such $\mu$ will have the form
\begin{equation}\label{p1skwformf}
\mu = \mu_0  + \alpha_1 (x \cdot t)n
 + \alpha_2 (x \cdot s)n 
 + \alpha_3 [(x \cdot t)s - (x \cdot s)t],
\end{equation}
where $\mu_0 \in \P_{0} \Lambda^0(f;\V)$.

Since $\P_{1,skw} \Lambda^0(f;\V)$ is a 6-dimensional space on each face,
the above quantities specify 24 degrees of freedom for the space $\P_{1,-}
\Lambda^2 (T;\V)$.  To see that these are a unisolvent set of degrees of
freedom for $\P_{1,-} \Lambda^2 (T;\V)$, we let $\omega = \omega_0 + d
\omega_1$, where $\omega_0 \in \P_{0}^+ \Lambda^2 (T;\V)$ and $\omega_1
\in \P_{1,f,-} \Lambda^1(T;\V)$ and set all degrees of freedom equal to
zero. Then for $\mu \in \P_{0} \Lambda^0(f;\V)$, since
\begin{equation*}
\int_f (\omega_0 + d \omega_1) \wedge \mu =\int_f  \omega_0 \wedge \mu,
\end{equation*}
we see that $\omega_0 =0$ by the unisolvency of the standard degrees of
freedom for $\P_{0}^+ \Lambda^2 (T;\V)$.  In addition, for
all $\mu \in \P_{1,skw} \Lambda^0(f;\V)$ and $\omega_0=0$,
\begin{equation*}
\int_f \omega \wedge \mu = \int_f d \omega_1 \wedge \mu
= \int_f \omega_1 \wedge d \mu.
\end{equation*}
Since $d \mu \in \P_0 \Lambda^1_{skw}(f;\V)$, $\omega_1 =0$ 
by the unisolvency of the degrees of freedom of the space
$\P_{1,f,-} \Lambda^1(T;\V)$.

Using an argument completely parallel to that used previously,
it is straightforward to show that the simplified spaces also
satisfy assumption \eqref{boxcond}, i.e., that $S_h$ is onto.

To translate the degrees of freedom of the space $\P_{1,-}\Lambda^2(T;
\V)$ to more standard finite element degrees of freedom, we use
the identification of an element 
$\omega \in \Lambda^2(T;\V)$ with a matrix  $F$ given by
$\omega(v_1,v_2) = F(v_1 \x v_2)$.  Then
$\omega(t,s) = F n$
and $\int_f \omega \wedge \mu = \int_f \mu^T F n \, df$.
Since $\mu \in \P_{1,skw} \Lambda^0(f;\V)$ and hence is of the form
\eqref{p1skwformf}, we get on each face the six degrees of freedom
\begin{equation*}
\int_f F n \,df, \quad 
\int_f (x \cdot t)n^T F n \, df, \quad
\int_f (x \cdot s)n^T F n \, df, \quad
\int_f [(x \cdot t)s^T  - (x \cdot s)t^T] F n \, df.
\end{equation*}

Finally, we note that the analogue of Theorem~\eqref{convergence-thm}
holds with $r=0$ for the simplified spaces.

\section*{Acknowledgments}
The authors are grateful to Geir Ellingsrud and Snorre
Christiansen at CMA and to Joachim Sch\"oberl, Johannes Kepler
University Linz, for many useful discussions.

\bibliographystyle{amsplain}
\bibliography{mixedelas3d-rev}

\providecommand{\bysame}{\leavevmode\hbox to3em{\hrulefill}\thinspace}
\begin{thebibliography}{10}

\bibitem{Adams-Cockburn}
Scot Adams and Bernardo Cockburn, \emph{A mixed finite element method
for elasticity in three dimensions}, Journal of Scientific Computing
\textbf{25} (2005), 515--521.

\bibitem{Amara-Thomas}
Mohamed Amara and Jean-Marie Thomas, \emph{Equilibrium finite elements for the
  linear elastic problem}, Numer. Math. \textbf{33} (1979), no.~4, 367--383.

\bibitem{Arnold-Awanou-Winther}
Douglas~N. Arnold, Gerard Awanou, and Ragnar Winther, \emph{Mixed
  finite elements for elasticity with strongly imposed symmetry}, in preparation.

\bibitem{Arnold-Brezzi-Douglas}
Douglas~N. Arnold, Franco Brezzi, and Jim Douglas, Jr., \emph{P{E}{E}{R}{S}: a
  new mixed finite element for plane elasticity}, Japan J. Appl. Math.
  \textbf{1} (1984), no.~2, 347--367.

\bibitem{Arnold-Douglas-Gupta}
Douglas~N. Arnold, Jim Douglas, Jr., and Chaitan~P. Gupta, \emph{A family of
  higher order mixed finite element methods for plane elasticity}, Numer. Math.
  \textbf{45} (1984), no.~1, 1--22.

\bibitem{Arnold-Falk}
Douglas~N. Arnold and Richard~S. Falk, \emph{A new mixed formulation for
  elasticity}, Numer. Math. \textbf{53} (1988), no.~1-2, 13--30.

\bibitem{Arnold-Falk-Winther05_1} Douglas~N. Arnold, Richard~S. Falk, and
Ragnar Winther \emph{Differential complexes and stability of finite element
methods. I: The de Rham complex}, in Compatible Spatial
Discretizations, D. Arnold, P. Bochev, R. Lehoucq, R. Nicolaides, and
M. Shashkov, eds., IMA Volumes in Mathematics and its Applications 142,
Springer Verlag 2005, 23-46.

\bibitem{Arnold-Falk-Winther05_2} Douglas~N. Arnold, Richard~S. Falk, and
Ragnar Winther \emph{Differential complexes and stability of finite element
methods. II: The elasticity complex}, in Compatible Spatial
Discretizations, D. Arnold, P. Bochev, R. Lehoucq, R. Nicolaides, and
M. Shashkov, eds., IMA Volumes in Mathematics and its Applications 142,
Springer Verlag 2005, 47-67.

\bibitem{acta} Douglas~N. Arnold, Richard~S. Falk, and
Ragnar Winther \emph{Finite element exterior calculus, homological techniques,
and application}, Acta Numerica (2006), 1-155.

\bibitem{Arnold-Winther02}
Douglas~N. Arnold and Ragnar Winther, \emph{Mixed finite elements for
  elasticity}, Numer. Math. \textbf{92} (2002), no.~3, 401--419.

\bibitem{BGG}
I.N. Bernstein, I.M. Gelfand, and S.I. Gelfand, \emph{Differential operators on
  the base affine space and a study of $\mathfrak g$--modules}, Lie groups and 
their representation, I.M. Gelfand, ed., (1975), 21--64.

\bibitem{Brezzi}
Franco Brezzi, \emph{On the existence, uniqueness and approximation of
  saddle-point problems arising from {L}agrangian multipliers}, Rev. Fran\c
  caise Automat. Informat. Recherche Op\'erationnelle S\'er. Rouge \textbf{8}
  (1974), no.~R-2, 129--151.

\bibitem{Brezzi-Fortin}
Franco Brezzi and Michel Fortin, \emph{Mixed and {H}ybrid {F}inite {E}lement
  {M}ethods}, Springer-Verlag, New York, 1991.

\bibitem{Cap}
Andreas \v{C}ap, Jan Slov\'{a}k, and Vladim\'{i}r Sou\v{c}ek,
\emph{Bernstein--Gelfand--Gelfand
sequences}, Ann. Math. (2) \textbf{154} (2001), 97--113.

\bibitem{Christensen}
Richard~M. Christensen, \emph{Theory of Viscoelasticity}, Dover Publications,
1982.

\bibitem{Ciarlet}
Philippe~G. Ciarlet, \emph{The finite element method for elliptic problems},
  North-Holland, Amsterdam, 1978.

\bibitem{douglas-roberts}
Jim Douglas, Jr. and Jean~E. Roberts, \emph{Global estimates for mixed 
methods for second order elliptic equations},  Math. Comp.  
\textbf{44}  (1985), no. 169, 39--52.

\bibitem{Eastwood}
Michael Eastwood, \emph{A complex from linear elasticity}, Rend. Circ. Mat.
Palermo (2) Suppl. (2000), no.~63, 23--29.

\bibitem{falk-cime1} Richard~S. Falk, \emph{Finite element methods for linear
    elasticity}, to appear in Mixed Finite Elements, Compatibility Conditions,
  and Applications, Lectures given at the C.I.M.E. Summer School held in
  Cetraro, Italy, June 26-July 1, 2006, Lecture Notes in Mathematics,
  Springer-Verlag.

\bibitem{falk-osborn}
Richard~S. Falk and John~E. Osborn, \emph{Error estimates for mixed 
methods},
R.A.I.R.O. Analyse num\'erique/Numerical Analysis, \textbf{14} (1980), no.~3,
 249--277.

\bibitem{Fraejisdv}
Baudoiun~M. Fraejis~de Veubeke, \emph{Stress function approach}, Proc.
of the World Congress on Finite Element Methods in Structural Mechanics, 
Vol. 1, Bournemouth, Dorset, England (Oct. 12-17, 1975), J.1--J.51.

\bibitem{Girault-Raviart}
V. Girault and P.-A. Raviart, \emph{Finite element methods for Navier-Stokes
equations. Theory and algorithms}, Springer Series in Computational
Mathematics, 5. Springer-Verlag, Berlin, 1986.

\bibitem{Johnson-Mercier}
Claes Johnson and Bertrand Mercier, \emph{Some equilibrium finite element
  methods for two-dimensional elasticity problems}, Numer. Math. \textbf{30}
  (1978), no.~1, 103--116.

\bibitem{Morley}
Mary~E. Morley, \emph{A family of mixed finite elements for linear elasticity},
  Numer. Math. \textbf{55} (1989), no.~6, 633--666.

\bibitem{Nedelec1}
Jean-Claude N\'ed\'elec, \emph{Mixed finite elements in $R^{3}$},
Numer. Math.  \textbf{35}  (1980), no. 3, 315--341. 

\bibitem{Nedelec2}
Jean-Claude N\'ed\'elec, \emph{A new family of mixed finite elements in
$R^3$}, Numer. Math. \textbf{50} (1986), no.~1, 57--81.

\bibitem{Stein-Rolfes}
Erwin Stein and Raimund Rolfes, \emph{Mechanical conditions for stability and
  optimal convergence of mixed finite elements for linear plane elasticity},
  Comput. Methods Appl. Mech. Engrg. \textbf{84} (1990), no.~1, 77--95.

\bibitem{Stenberg86}
Rolf Stenberg, \emph{On the construction of optimal mixed finite element
  methods for the linear elasticity problem}, Numer. Math. \textbf{48} (1986),
  no.~4, 447--462.

\bibitem{Stenberg88}
\bysame, \emph{A family of mixed finite elements for the elasticity problem},
  Numer. Math. \textbf{53} (1988), no.~5, 513--538.

\bibitem{Stenberg-mafelap}
\bysame, \emph{Two low-order mixed methods for the elasticity problem}, The
  mathematics of finite elements and applications, VI (Uxbridge, 1987),
  Academic Press, London, 1988, pp.~271--280.

\bibitem{Watwood-Hartz}
Vernon~B. Watwood~Jr. and B.~J. Hartz, \emph{An equilibrium stress field model
  for finite element solution of two--dimensional elastostatic problems},
  Internat. Jour. Solids and Structures \textbf{4} (1968), 857--873.

\end{thebibliography}

\end{document}